# An Improved Milstein Method for the Numerical Solution of Multidimensional Stochastic Differential Equations

Paromita Banerjee[*] & Anirban Mondal[†]

January 11, 2026

**Abstract**

Stochastic differential equations (SDEs) offer powerful and accessible mathematical models for capturing both deterministic and probabilistic aspects of dynamic behavior across a wide range of physical, financial, and social systems. However, analytical solutions for many SDEs are often unavailable, necessitating the use of numerical approximation methods. The rate of convergence of such numerical methods is of great importance, as it directly influences both computational efficiency and accuracy. This paper presents a proposed theorem, along with its proof, that facilitates the numerical evaluation of the strong (and weak) order of convergence of a numerical scheme for an SDE when the analytical solution is unavailable. Additionally, we address the challenge of numerically computing the multiple stochastic integrals required by the Milstein method to achieve improved convergence rates for multidimensional SDEs. In this context, two newly proposed numerical techniques for computing these multiple stochastic integrals are introduced and compared with existing approaches in terms of efficiency and effectiveness. The methodologies are further illustrated through simulation studies and applications to widely used financial models.

## 1 Introduction

Approximately three centuries ago, Newton and Leibniz developed differential calculus, providing a mathematical framework for modeling continuous-time dynamical systems in mechanics, astronomy, and numerous other scientific disciplines. This foundational calculus has supported the transformative advances in science, technology, and manufacturing witnessed over the last two centuries. As efforts have progressed toward constructing more realistic models, it has become essential to incorporate stochastic effects. In fields such as finance, for example, randomness is not merely a refinement but a central feature of the system dynamics. Continuous-time stochastic models now play a critical role in diverse areas of application, including microelectronics, signal processing and filtering, various branches of biology and physics, population dynamics, epidemiology, psychology, economics, finance, insurance, fluid dynamics, radio astronomy, hydrology, structural mechanics, chemistry, and medicine. Practical problems arising in several of these domains during the mid-20th century motivated the development of a corresponding stochastic calculus. Following Bachelier's early use of Brownian motion in modeling stock prices in the Paris Bourse [2] and Einstein's mathematically equivalent formulation, Wiener [19] developed a more complete mathematical theory of Brownian motion. A major advancement came with Itô's seminal work [8], which established the foundations of what is now known as Itô calculus. This stochastic extension of classical differential calculus enables continuous-time modeling of phenomena such as stock price dynamics or the motion of microscopic particles under random perturbations. The resulting stochastic differential equations (SDEs) generalize ordinary differential equations (ODEs) by incorporating random forcing terms. SDEs driven by Brownian motion or Lévy processes now serve as essential tools across a wide spectrum of applications, including chemistry, biology, epidemiology, mechanics, microelectronics, economics, and finance.

We begin with the basic one-dimensional stochastic differential equation (SDE). Throughout this paper, we considered a filtered probability space $(\Omega, \mathscr{F}, \mathbb{F} = \{\mathscr{F}_t\}_{t \geq 0}, P)$ satisfying the usual hypotheses (right-continuity and completeness). A one-dimensional Itô SDE driven by a Wiener process, $W_t$, takes the form

$$dX_t = a(t, X_t)\,dt + b(t, X_t)\,dW_t, \qquad t \in [t_0, T], \tag{1}$$

where the initial value $X_{t_0}$ is $\mathscr{F}_{t_0}$-measurable, independent of $\{W_t\}_{t \geq 0}$, and satisfies $E(X_{t_0}^2) < \infty$. Since Brownian motion is continuous but nowhere differentiable, (1) cannot be interpreted in the classical derivative sense and is instead defined through the equivalent integral equation

$$X_t = X_{t_0} + \int_{t_0}^{t} a(s, X_s)\,ds + \int_{t_0}^{t} b(s, X_s)\,dW_s. \tag{2}$$

The first integral is an ordinary Riemann integral, whereas the second is the Itô stochastic integral. For a partition $c = t_0 < t_1 < \cdots < t_n = d$ of $[c,d]$, the Riemann integral is defined as $\int_c^d f(x)\,dx = \lim_{\Delta t \to 0} \sum_{i=1}^{n} f(t_i') \Delta t_i$, with $t_{i-1} \leq t_i' \leq t_i$ and $\Delta t_i = t_i - t_{i-1}$.

---

[*]Department of Mathemtics, Computer Science, and Data Science, John Carroll University
[†]Department of Mathematics, Applied Mathematics, and Statistics, Case Western Reserve University

In contrast, the Itô integral is defined using left-endpoint evaluations: $\int_c^d f(t)\,dW_t = \lim_{\Delta t \to 0} \sum_{i=1}^n f(t_{i-1})\Delta W_i$, where $\Delta W_i = W_{t_i} - W_{t_{i-1}}$. Because both $f(t)$ and $W_t$ are random, the integral is itself a random variable. The differential notation is a standard shorthand: $dI = f(t)\,dW_t \iff I = \int_c^d f(t)\,dW_t$. The term $dW_t$ is often interpreted as idealized white noise, and solutions of (1) naturally combine deterministic drift with stochastic diffusion generated by Brownian motion.

**Existence and uniqueness** ([9], Thm. 5.4): If $a$ and $b$ satisfy either the global or local Lipschitz condition together with the linear growth condition, and if $X_{t_0}$ is $\mathscr{F}_{t_0}$-measurable, independent of $W$, with $E[\|X_{t_0}\|^2] < \infty$, then the SDE 1 has a unique strong solution $\{X_t\}_{t_0 \le t \le T}$ satisfying $E\left[\sup_{t_0 \le t \le T} \|X_t\|^2\right] \le C(1 + E[\|X_{t_0}\|^2])$, where $C > 0$ depends only on $K$ and $T$.

Analytical solutions of SDEs rely on the stochastic chain rule, known as *Itô's formula*, introduced by Itô [8]. Let $g$ be a nonanticipative, mean-square integrable process, so that $g_{t_n}$ is independent of $W_{t_{n+1}} - W_{t_n}$. Then the Itô integral satisfies $E\left(\int_{t_0}^T g_s\,dW_s\right) = 0, \quad E\left(\int_{t_0}^T g_s\,dW_s\right)^2 = \int_{t_0}^T E[g_s^2]\,ds$.

**Itô's Lemma.** For the scalar SDE $dX_t = a(X_t)\,dt + b(X_t)\,dW_t$ and a function $f(t,x)$ with continuous $\partial_t f$ and $\partial_x^2 f$, we have

$$df(t, X_t) = \left(\partial_t f + a(X_t)\partial_x f + \tfrac{1}{2} b^2(X_t)\partial_x^2 f\right)dt + b(X_t)\partial_x f\,dW_t. \tag{3}$$

The additional drift term $\tfrac{1}{2}b^2 \partial_x^2 f$ is characteristic of Itô calculus and plays a central role in the construction of numerical schemes for SDEs. In integral form,

$$f(t, X_t) = f(t_0, X_{t_0}) + \int_{t_0}^t L^0 f(s, X_s)\,ds + \int_{t_0}^t L^1 f(s, X_s)\,dW_s, \tag{4}$$

where $L^0 = \partial_t + a\partial_x + \tfrac{1}{2}b^2 \partial_x^2$ and $L^1 = b\partial_x$.

An alternative stochastic integral is the *Stratonovich integral* [17], where the integrand is evaluated at the midpoint of each interval. Unlike Itô integrals, Stratonovich integrals obey the standard chain rule, allowing deterministic calculus techniques. Both calculi are mathematically valid, and a suitable drift modification can convert one to the other. In this work, we focus on Itô calculus due to its convenient expectation properties and widespread use in numerical SDE analysis.

**Multivariate SDE:** Next, we consider a multivariate SDE of the following form

$$d\mathbf{X}_t = a(t, \mathbf{X}_t)\,dt + b(t, \mathbf{X}_t)\,d\mathbf{W}_t, \quad \mathbf{X}_0 = x_0. \tag{5}$$

where, $\mathbf{X}_t = (X_t^1, \ldots, X_t^d)$, $\mathbf{W}_t$ is a $m$-dimensional Wiener process, $a : \mathbb{R}_+ \times \mathbb{R}^d \to \mathbb{R}^d$, is a Drift function, $b : \mathbb{R}_+ \times \mathbb{R}^d \to M(d, m)$, is a Diffusion function, $M(d,m)$ denote the set of $d \times m$ real matrices, and $x_0 \in \mathbb{R}^d$, is the Initial vector. The process $\mathbf{X}_t = \mathbf{X}(t)$ is called a strong solution of the SDE 5 if for all $t > 0$, $\mathbf{X}_t$ is a function $F(t, (\mathbf{W}_s, s \le t))$ of the given Wiener process $\mathbf{W}$, the integrals $\int_0^t a(s, \mathbf{X}_s)\,ds$ and $\int_0^t b(s, \mathbf{X}_s)\,d\mathbf{W}_s$ exist, and

$$\mathbf{X}_t = x_0 + \int_0^t a(s, \mathbf{X}_s)\,ds + \int_0^t b(s, \mathbf{X}_s)\,d\mathbf{W}_s. \tag{6}$$

Assume the coefficients $a$ and $b$ satisfy the Lipschitz condition and the linear growth condition, and $x_0$ is constant, then there exists a unique strong solution $\mathbf{X}_t$ with continuous trajectories satisfying $E[\|\mathbf{X}_t\|^2] \le Ce^{Ct}(1 + \|x_0\|^2)$, where $C > 0$ depends on $K$ and $T$.

**Multivariate Itô lemma:** Let $\mathbf{X}_t$ satisfy

$$dX_t^i = a_i(X_t)\,dt + \sum_{j=1}^m b^{ij}(X_t)\,dW_t^j, \quad i = 1, \ldots, d, \tag{7}$$

with $b^{ij}(x)$'s are such that the covariance matrix $\Sigma = [\sum_{l=1}^m b^{il} b^{jl}]$ is positive definite. Then for $f : \mathbb{R}_+ \times \mathbb{R}^d \to \mathbb{R}$ with sufficient differentiability,

$$df(t, X_t) = \left(\frac{\partial f}{\partial t} + \sum_{i=1}^d a_i(X_t)\frac{\partial f}{\partial x_i} + \frac{1}{2}\sum_{i,j=1}^d \sum_{l=1}^m b^{il}(X_t)b^{jl}(X_t)\frac{\partial^2 f}{\partial x_i \partial x_j}\right)dt + \sum_{i=1}^d \sum_{l=1}^m b^{il}(X_t)\frac{\partial f}{\partial x_i}\,dW_t^l. \tag{8}$$

For correlated Brownian motions, additional correlation terms appear in the diffusion part of the Itô formula. Suppose there is correlation matrix is given by $\rho = ((\rho_{ij}))$, where $\rho_{ij} = \rho_{ji}$ are such that $\frac{1}{t}\mathbb{E}(W_t^i W_t^j) = \rho_{ij}\ \forall i, j \le n$. Then equation 8 becomes

$$df(t, X_t) = \left(\frac{\partial}{\partial t} + \sum_{i=1}^m a^i(X_t)\frac{\partial}{\partial x_i} + \frac{1}{2}\sum_{i,j=1}^m \sum_{l,k=1}^n \rho_{lk} b^{il}(X_t) b^{kj}(X_t)\frac{\partial^2}{\partial x_i \partial x_j}\right)f(t, X_t)\,dt + \sum_{i=1}^m \sum_{l=1}^n b^{il}(X_t)\frac{\partial}{\partial x_i}f(t, X_t)\,dW_t^l \tag{9}$$

The existence and uniqueness conditions still hold, provided $a_i$ and $b^{ij}$ satisfy Lipschitz and linear growth assumptions.

In many cases, analytic solutions of such SDEs are unavailable, necessitating numerical approximations. These approximations are typically based on a *time discretization* of the interval $[0, T]$: $0 = t_0 < t_1 < \cdots < t_n < \cdots < t_N = T$, with step-size $\Delta_n = t_{n+1} - t_n$. More general time discretizations, including random steps, can also be considered. Simulation experiments and theoretical studies have shown that not all classical or heuristic time-discrete approximations converge to the true solution as $\Delta_n \to 0$ [5, 20, 6, 4, 16]. This highlights the need for systematic selection of efficient and reliable numerical methods. Among discrete-time approximations, the most widely used methods are the Euler-Maruyama [14] method and the Milstein method [15].

When approximating SDEs numerically, it is important to quantify how well the discrete-time solution $\{X_n\}_{n=0}^{N}$ approximates the true continuous-time solution $\{X_t\}_{t\in[0,T]}$. Two commonly used notions of convergence are *strong convergence* and *weak convergence*: Strong convergence measures the accuracy of the sample paths of the numerical solution. Weak convergence captures the accuracy of the distribution of the numerical solution rather than individual paths.

In practice, the Euler-Maruyama method typically achieves strong order 0.5 and weak order 1, while the Milstein method improves the strong order to 1 for scalar SDEs under sufficient smoothness conditions [11]. These convergence rates play a crucial role in the efficiency and reliability of the numerical scheme for a given SDE.

For a given SDE, the strong and weak rates of convergence of a numerical scheme are often studied numerically via computer simulation experiments. However, such simulation-based approaches typically require the true analytical solution to compute the expected value of the absolute error, which is generally unavailable. To overcome this limitation, we present and prove theorems that enable the numerical computation of strong and weak convergence rates without requiring the analytical solution of the SDEs. The theorems emphasize the relationship between the time discretization stepsize and the order of convergence (strong, weak, and mean-square) of the numerical method.

The Euler–Maruyama numerical scheme for a scalar SDE extends naturally to multivariate SDEs by applying the one-dimensional update to each component separately. However, as in the scalar case, the Euler–Maruyama scheme for multivariate SDEs attains only strong order 0.5. To achieve strong order 1, one must use the Milstein scheme for multivariate SDEs. The main difficulty in implementing the Milstein method in higher dimensions lies in the computation of the Itô multiple stochastic integrals, which cannot be expressed in terms of the increments of the Wiener process components, nor is there a known tractable probability distribution for these integrals. In this paper, we examine in detail several methodologies for approximating these multiple stochastic integrals, together with the resulting convergence properties of the corresponding Milstein schemes. [13] introduced the Brownian bridge process associated with an $m$-dimensional Wiener process, and [11] proposed approximating the multiple stochastic integral by first expressing each component of the Wiener process via a Fourier expansion of the Brownian bridge, and then integrating iteratively with respect to each component. However, this approach requires generating an additional $2(2p+1)$ normal random variables, which leads to high variability in the estimator and significant computational cost for small $\Delta$. Moreover, because the method involves two layers of approximation, its validity is restricted to sufficiently small step-sizes $\Delta$, limiting its practical applicability. [10] suggested an alternative computational strategy in which the multiple stochastic integral inside the Milstein scheme is approximated using a subdivision (or "inner-loop") method that solves another SDE using a two-dimensional Euler–Maruyama method. However, [4] showed that, with the initial conditions as specified in [10], the resulting inner-loop solution does not yield the desired stochastic process. We show that modifying the initial conditions to zero in each dimension leads to the correct approximation. We refer to this corrected method as "*E-M IC= 0*". We also derive the sufficient condition for strong order 1 convergence of the parent Milstein scheme and present simulation results supporting our theoretical findings. Furthermore, we propose a new method for approximating the multiple stochastic integral, inside the Milstein method, using a related subdivision-based approach, in which the inner-loop SDE is solved using an approximate version of the Milstein method itself. The proposed scheme yields an improved mean-square error convergence rate compared with the "*E-M IC= 0*" method, and achieves the same strong order of convergence for the parent Milstein scheme while requiring only half as many subintervals. These results are established both theoretically and through numerical experiments.

We apply the proposed Milstein schemes to several SDEs commonly used in financial modeling. For the one-dimensional case, we consider the classical Black–Scholes model [3]. For the two-dimensional setting, we apply and compare the proposed Milstein schemes on the generalized Heston stochastic volatility model [7, 1].

The paper is organized as follows. Section 2 presents the rates of convergence of numerical schemes for solving SDEs, along with the proposed theorems and their proof for evaluating these convergence rates. Section 3 focuses on numerical methods for multidimensional SDEs, comparing existing approaches with the proposed variants of the multidimensional Milstein scheme. Section 4 provides applications to SDEs in financial models. Finally, Section 5 concludes the paper with a brief discussion and potential directions for future research.

## 2 Numerical Schemes for SDEs and Rates of Convergence

We begin by describing the most commonly used numerical scheme for one-dimensional SDEs - the Euler-Maruyama method.

### 2.1 Euler–Maruyama Method

The Euler approximation, first studied in [14], is the simplest numerical method for SDEs and is well suited for implementation on digital computers. The stochastic analogue of the Euler scheme for the Itô SDE (1) is given by the recursion

$$X_{n+1} = X_n + a(t_n, X_n)\Delta_n + b(t_n, X_n)\Delta W_n, \tag{10}$$

where $\Delta_n = t_{n+1} - t_n \int_{t_n}^{t_{n+1}} ds = I_{(0)}$ is the is the length of the time discretization interval $[t_n, t_{n+1}]$ for $n = 0, 1, \ldots, N-1$, with $X_0 = X_{t_0}$ and $X_n = X_{t_n}$. The noise increments are defined as $\Delta W_n = W_{t_{n+1}} - W_{t_n}$, and are independent Gaussian random variables with distribution $N(0, \Delta_n)$.

It is well known that the Euler approximation converges to the true solution $X_t$ of the Itô SDE (1) as $\Delta_n \to 0$, under various notions of convergence. To simulate the Euler–Maruyama method, one must generate the independent Gaussian increments $\Delta W_n$, typically using pseudo-random number generators. As in deterministic numerical analysis, the Euler method is simple and computationally inexpensive, but may be inefficient and can exhibit poor stability properties. The Euler–Maruyama scheme is consistent with the Itô calculus because the stochastic term in (3) approximates the Itô integral in (2) over the interval $[t_n, t_{n+1}]$ by evaluating its integrand at the left endpoint: $\int_{t_n}^{t_{n+1}} b(s, X_s)\, dW_s \approx \int_{t_n}^{t_{n+1}} b(t_n, X_n)\, dW_s = b(t_n, X_n) \int_{t_n}^{t_{n+1}} dW_s$.

## 2.2 Convergence of Numerical Methods for SDEs

It is useful to classify the efficiency of a numerical scheme by identifying its order of convergence. Convergence for numerical schemes for SDEs can be defined in several meaningful ways. In particular, it is common to distinguish between *strong* and *weak* convergence, depending on whether one requires the realizations of the paths or only their probability distributions to be close. Unlike the deterministic setting, the stochastic environment admits several types of convergence that are theoretically or practically relevant. Therefore, in stochastic numerical analysis, one must specify the class of problems of interest before designing a numerical method and determining which convergence criterion is most suitable. In practice, two principal types of convergence are typically considered: i) approximation of *sample paths*, and ii) approximation of the *probability distribution*. For the purposes of classifying numerical algorithms, we adopt simplified characterizations of these two types of convergence, referred to as the strong and weak convergence criteria, respectively.

### 2.2.1 Strong Order of Convergence

Tasks that involve direct simulation of sample paths, such as generating stock price scenarios in finance, computing filtering estimates for hidden variables, or evaluating statistical estimators for SDE parameters, require that the simulated sample paths closely approximate those of the true solution. In such applications, a strong convergence criterion is appropriate.

Under suitable assumptions on the SDE, consider a fixed interval $[t_0, T]$ and let $\Delta$ denote the maximum step size of a discretization of this interval. A numerical scheme has *strong order of convergence* $\gamma \in (0, \infty]$ if there exists a constant $K < \infty$ such that

$$E\left|X_T - X_T^\Delta\right| \leq K\Delta^\gamma, \tag{11}$$

where $X_T$ is the exact solution and $X_T^\Delta$ is the numerical approximation.

The Euler-Maruyama (EM) scheme corresponds to truncating the Itô–Taylor expansion after including terms involving only time increments and Wiener increments of multiplicity one. Thus, the Euler approximation leads to a strong Itô–Taylor approximation of order $\gamma = \frac{1}{2}$. We next state standard conditions ensuring strong convergence of the Euler scheme with order $\gamma = \frac{1}{2}$.

**Euler Scheme (Strong Convergence).** Assume $E|X_0|^2 < \infty$ and $E\left(|X_0 - X_0^\Delta|^2\right)^{1/2} \leq K_1 \Delta^{1/2}$. Suppose the drift $a$ and diffusion $b$ satisfy, for all $s, t \in [0, T]$ and $x, y \in \mathbb{R}^d$, $|a(t, x) - a(t, y)| + |b(t, x) - b(t, y)| \leq K_2|x - y|$, $|a(t, x)| + |b(t, x)| \leq K_3(1 + |x|)$, $|a(s, x) - a(t, x)| + |b(s, x) - b(t, x)| \leq K_4(1 + |x|)|s - t|^{1/2}$, where $K_1, \ldots, K_4$ are independent of $\Delta$. Then the Euler approximation $X_T^\Delta$ satisfies $E\left(|X_T - X_T^\Delta|^2\right)^{1/2} \leq K_5 \Delta^{1/2}$, with $K_5$ independent of $\Delta$. A proof is given in [12], Section 10.2.2.

### 2.2.2 Weak Order of Convergence

In many applications, a pathwise approximation of the solution is unnecessary. For example, when computing moments of $X_T$, probabilities associated with $X_T$, option prices, or general functionals of the form $E[g(X_T)]$, only the distribution of $X_T$ needs to be approximated accurately. In such cases, a much weaker form of convergence suffices.

A numerical approximation $X_T^\Delta$ of a solution $X_T$ of an SDE is said to converge *weakly* with order $\beta \in (0, \infty]$ if, for any polynomial function $g$, there exists a constant $K_g < \infty$ such that

$$\left|E[g(X_T)] - E[g(X_T^\Delta)]\right| \leq K_g \Delta^\beta, \tag{12}$$

whenever the expectations exist.

The strong and weak orders of convergence correspond to the largest possible values of $\gamma$ and $\beta$ in (11) and (12), respectively. This criterion includes the convergence of $p$th moments by taking $g(x) = x^p$, and reduces to the deterministic convergence criterion in the case $b = 0$ and $g(x) = x$.

In this work, we primarily focus on the strong rate of convergence of the Euler–Maruyama and Milstein schemes. The Euler-Maruyama scheme has strong order $\gamma = \frac{1}{2}$ and weak order $\beta = 1$ for SDEs with sufficiently smooth coefficients (e.g., continuously differentiable with uniformly bounded derivatives). Higher orders may be possible for restricted classes of SDEs, such as those with additive noise (where $b$ does not depend on the state $x$), which can achieve strong order $\gamma = 1$. However, the strong order $\gamma = \frac{1}{2}$ and weak order $\beta = 1$ of the Euler scheme are relatively low, especially given that many applications require generating a large number of sample paths. This motivated the development and study of higher-order numerical schemes.

## 2.3 Milstein Method

A higher order of convergence cannot be obtained using deterministic numerical schemes for SDEs, even when they are consistent, because such schemes involve only the simple increments of time $\Delta_n$ and noise $\Delta W_n$. The latter provides a poor approximation of the highly irregular Wiener process within the discretization subinterval $[t_n, t_{n+1}]$. Achieving higher-order convergence requires incorporating additional information about the Wiener process over each subinterval. This information is supplied by multiple stochastic integrals, which arise in stochastic Taylor expansions of SDE solutions. Numerical schemes of arbitrarily high order can be derived by truncating suitable stochastic Taylor expansions, themselves obtained through iterated application of the Itô formula (4).

When $f(t, x) = x$, the Itô formula (4) reduces to the integral form (2) of the SDE (1), i.e.,

$$X_t = X_{t_0} + \int_{t_0}^t a(s, X_s)\,ds + \int_{t_0}^t b(s, X_s)\,dW_s. \tag{13}$$

Applying the Itô formula to the integrand functions $f(t,x) = a(t,x)$ and $f(t,x) = b(t,x)$ in (13) yields

$$X_t = X_{t_0} + \int_{t_0}^t \left[ a(t_0, X_{t_0}) + \int_{t_0}^s L^0 a(u, X_u)\,du + \int_{t_0}^s L^1 a(u, X_u)\,dW_u \right] ds + \int_{t_0}^t \left[ b(t_0, X_{t_0}) + \int_{t_0}^s L^0 b(u, X_u)\,du + \int_{t_0}^s L^1 b(u, X_u)\,dW_u \right] dW_s$$

$$= X_{t_0} + a(t_0, X_{t_0}) \int_{t_0}^t ds + b(t_0, X_{t_0}) \int_{t_0}^t dW_s + R_1(t, t_0). \quad (14)$$

The remainder term is

$$R_1(t, t_0) = \int_{t_0}^s \int_{t_0}^s L^0 a(u, X_u)\,du\,ds + \int_{t_0}^s \int_{t_0}^s L^1 a(u, X_u)\,dW_u\,ds + \int_{t_0}^t \int_{t_0}^s L^0 b(u, X_u)\,du\,dW_s + \int_{t_0}^t \int_{t_0}^s L^1 b(u, X_u)\,dW_u\,dW_s. \quad (15)$$

Replacing $t_0$ by $t_n$, $X_{t_0}$ by $X_n$, and $t$ by $t_{n+1}$, then discarding the remainder, yields the Euler scheme (10), the simplest non-trivial stochastic Taylor scheme.

Higher-order stochastic Taylor expansions are obtained by repeatedly applying the Itô formula to the integrand functions within the remainder terms. For instance, applying the Itô formula to the integrand $L^1 b$ in the fourth double integral of $R_1(t, t_0)$ produces

$$X_t = X_{t_0} + a(t_0, X_{t_0}) \int_{t_0}^t ds + b(t_0, X_{t_0}) \int_{t_0}^t dW_s + L^1 b(t_0, X_{t_0}) \int_{t_0}^t \int_{t_0}^s dW_u\,dW_s + R_2(t, t_0), \quad (16)$$

where,

$$R_2(t, t_0) = \int_{t_0}^s \int_{t_0}^s L^0 a(u, X_u)\,du\,ds + \int_{t_0}^s \int_{t_0}^s L^1 a(u, X_u)\,dW_u\,ds + \int_{t_0}^t \int_{t_0}^s L^0 b(u, X_u)\,du\,dW_s + \int_{t_0}^t \int_{t_0}^s \int_{t_0}^u L^0 L^1 b(v, X_v)\,dv\,dW_u\,dW_s$$

$$+ \int_{t_0}^t \int_{t_0}^s \int_{t_0}^u L^1 L^1 b(v, X_v)\,dW_v\,dW_u\,dW_s. \quad (17)$$

Substituting $t_0 = t_n$, $X_{t_n} = X_n$, and $t = t_{n+1}$, and discarding the remainder, yields the *Milstein scheme*

$$X_{n+1} = X_n + a(t_n, X_n)\Delta_n + b(t_n, X_n)\Delta W_n + L^1 b(t_n, X_n) \int_{t_n}^{t_{n+1}} \int_{t_n}^s dW_u\,dW_s, \quad (18)$$

originally proposed by Milstein in the 1970s [15]. The Milstein method has strong order $\gamma = 1$ and weak order $\beta = 1$. Thus, it improves the strong convergence order compared to Euler–Maruyama but does not improve the weak order. To obtain higher-order schemes, this procedure can be continued by expanding further integrand functions using the Itô formula. As suggested by the Euler and Milstein schemes: i) higher-order schemes obtain their improved order by incorporating multiple stochastic integrals, and ii) a scheme may possess different strong and weak convergence orders.

A compact multi-index notation introduced by Wagner [18] clarifies which stochastic integrals must be included to obtain a desired order. Denoting $b^0(t,x) = a(t,x)$, $dW_t^0 = dt$, $b^1(t,x) = b(t,x)$, and $dW_t^1 = dW_t$, the Itô formula becomes

$$f(t, X_t) = f(t_0, X_{t_0}) + \sum_{j=0}^1 \int_{t_0}^t L^j f(s, X_s)\,dW_s^j,$$

and for $f(t,x) = x$ yields the SDE

$$X_t = X_{t_0} + \sum_{j=0}^1 b^j(t, X_t)\,dW_t^j.$$

In this notation, the coefficients of the Euler scheme correspond to $L^j \mathrm{id} = b^j$ with integrals $\int_{t_n}^{t_{n+1}} dW_s^j$, for $j = 0, 1$. The Milstein scheme additionally includes $L^1 L^1 \mathrm{id} = L^1 b^1$ with associated integral $\int_{t_n}^{t_{n+1}} \int_{t_n}^s dW_u^1\,dW_s^1$, corresponding to the multi-index $(1, 1)$.

Approximating multiple stochastic integrals of higher order requires additional care in simulation. Fortunately, some have exact formulas, such as

$$I_{j,j}[t_n, t_{n+1}] = \int_{t_n}^{t_{n+1}} \int_{t_n}^t dW_s^j\,dW_t^j = \frac{1}{2}\{(\Delta W_n^j)^2 - \Delta_n\}. \quad (19)$$

Others have known distributions. For example, $\Delta Z_n^j = I_{j,0}[t_n, t_{n+1}] = \int_{t_n}^{t_{n+1}} \int_{t_n}^t dW_s^j\,dt$, is Gaussian with mean zero and variance $\frac{1}{3}\Delta_n^3$, and $E[\Delta W_n^j \Delta Z_n^j] = \frac{1}{2}\Delta_n^2$, for all $j \geq 1$.

Using (19), the one-dimensional Milstein scheme becomes

$$X_{n+1} = X_n + a\Delta_n + b\Delta W_n + \frac{1}{2} b b' \left( (\Delta W_n)^2 - \Delta_n \right). \quad (20)$$

**Strong convergence of Milstein scheme:** Assume $E|X_0|^2 < \infty$ and $E|X_0 - X_0^\Delta|^{2^{1/2}} \leq K_1 \Delta^{1/2}$. For all $t, s \in [0, T]$, $x \in \mathbb{R}^m$, $j = 0, \ldots, d$, $j_1, j_2 = 1, \ldots, d$, let the coefficients satisfy the Lipschitz and linear growth conditions $|a(t,x) - a(t, x^\Delta)| + |b^{j_1}(t,x) - b^{j_1}(t, x^\Delta)| + |L^{j_1} b^{j_2}(t,x) - L^{j_1} b^{j_2}(t, x^\Delta)| \leq K_2 |x - x^\Delta|$, $|a(t,x)| + |L^j a(t,x)| + |b^{j_1}(t,x)| + |L^j b^{j_2}(t,x)| + |L^{j_1} L^{j_1} b^{j_2}(t,x)| \leq K_3(1 + |x|)$, and the temporal Hölder conditions $|a(s,x) - a(t,x)| + |b^{j_1}(s,x) - b^{j_1}(t,x)| + |L^{j_1} b^{j_2}(s,x) - L^{j_1} b^{j_2}(t,x)| \leq K_4(1 + |x|)|s - t|^{1/2}$, where $K_1, \ldots, K_4$ are independent of $\Delta$. Then the Milstein approximation satisfies $E(|X_T - X_T^\Delta|) \leq K_5 \Delta$, with $K_5$ independent of $\Delta$. The proof follows the standard Itô–Taylor expansion argument (see [11]), and the result implies first-order strong convergence of the Milstein scheme.

## 2.4 Numerical Study of the Rate of Convergence with Unknown True Solution

For a given SDE, the strong rate of convergence of a numerical scheme, as defined in (11), can be studied numerically through simulation. Suppose we aim to solve (2) with initial value $X_0$. To estimate the strong rate of convergence, we can proceed through the following steps:

1. Fix a simulation interval $[0, T]$ and choose a time step $\Delta = T/N$, where $N$ is the number of time steps. Set the initial value $X_0$.
2. Generate increments of a Wiener process, $\Delta W_n \overset{\text{i.i.d.}}{\sim} N(0, \Delta)$, for $n = 1, \ldots, N$.
3. Using a numerical scheme (e.g., Euler–Maruyama), update the solution by
$$X_{n+1} = X_n + a(t_n, X_n)\Delta + b(t_n, X_n)\Delta W_n,$$
   and obtain $X_T^\Delta = X_N$. For other schemes, use the corresponding update formula.
4. If the analytical solution of (2) is known, compute the exact value $X_T$ using the *same* Wiener increments.
5. Compute the absolute error $\text{Error}(\Delta) = |X_T - X_T^\Delta|$.
6. Repeat Steps 2–5 for $M$ independent realizations of the Wiener process to obtain errors $\text{Error}_i(\Delta)$, $i = 1, \ldots, M$.
7. Estimate the expected value of the absolute error by the Monte Carlo average $\widehat{E}(\Delta) = \frac{1}{M} \sum_{i=1}^{M} \text{Error}_i(\Delta)$.

The above steps yield an empirical estimate of $E|X_T - X_T^\Delta|$.

Now from (11), we have $\log(E|X_T - X_T^\Delta|) \leq \log(K) + \gamma \log(\Delta)$, which motivates the following algorithm for estimating the strong convergence rate $\gamma$.

1. Define a sequence of time steps $\Delta_r = T/N_r$, $r = 1, \ldots, R$.
2. For each $\Delta_r$, compute the Monte Carlo estimate $\widehat{E}(\Delta_r)$ using the steps described above. If the time steps satisfy $N_r = rN$ for a fixed $N$, then Wiener increments for coarser grids can be constructed from those of the finest grid: $\Delta_r W_n = \sum_{j=(n-1)k+1}^{nk} \Delta_R W_j$, $k = N_R/N_r$.
3. Fit a linear regression model with $\log(\widehat{E}(\Delta_r))$ as the response and $\log(\Delta_r)$ as the covariate. The slope of the fitted line provides an estimate of the strong convergence rate $\gamma$.

In principle, this simulation approach requires knowledge of the true solution $X_T$ of the SDE. However, for many SDEs the analytical solution is not available, which is often the reason for employing numerical schemes in the first place. To address this challenge, we employ the following result.

**Theorem 2.1.** *Let $X_T^\Delta$ denote a time-discrete approximation with step size $\Delta$ of the solution $X_T$ of an SDE, constructed using a given Wiener process. Then*
$$E\left|X_T - X_T^\Delta\right| < C_1 \Delta^\gamma, \tag{21}$$
*for all sufficiently small $\Delta$, with $C_1$ independent of $\Delta$, if and only if the scheme is strongly convergent and, for any integer $k > 1$, there exists a constant $C_2$ independent of $\Delta$ such that*
$$E\left|X_T^\Delta - X_T^{\Delta/k}\right| < C_2 \Delta^\gamma, \tag{22}$$
*where $X_T^{\Delta/k}$ denotes the approximation using the finer step size $\Delta/k$ and the same Wiener process.*

*Proof.* Assume that (21) holds for all $\Delta > 0$. Then, for any $k > 1$, we have
$$E\left|X_T - X_T^{\Delta/k}\right| < C_1 \frac{\Delta^\gamma}{k^\gamma}. \tag{23}$$

Now,
$$E\left|X_T^\Delta - X_T^{\Delta/k}\right| = E\left|(X_T - X_T^{\Delta/k}) - (X_T - X_T^\Delta)\right| \leq E\left|X_T - X_T^{\Delta/k}\right| + E\left|X_T - X_T^\Delta\right| < C_1 \frac{\Delta^\gamma}{k^\gamma} + C_1 \Delta^\gamma$$
$$= C_2 \Delta^\gamma,$$

where $C_2 = C_1\left(1 + \frac{1}{k^\gamma}\right)$, which is independent of $\Delta$. Moreover,
$$E\left|X_T - X_T^\Delta\right| < C_1 \Delta^\gamma \implies \lim_{\Delta \to 0} E\left|X_T - X_T^\Delta\right| = 0. \tag{24}$$

This proves one direction of the theorem.

Now assume that (22) holds for all $\Delta > 0$ and that $X_T^\Delta$ is strongly convergent. We write
$$X_T - X_T^\Delta = X_T - X_T^{\Delta/k^M} + \sum_{i=0}^{M-1} \left(X_T^{\Delta/k^{i+1}} - X_T^{\Delta/k^i}\right).$$

Hence,
$$E|X_T - X_T^\Delta| \le E\left|X_T - X_T^{\Delta/k^M}\right| + \sum_{i=0}^{M-1} E\left|X_T^{\Delta/k^{i+1}} - X_T^{\Delta/k^i}\right|.$$

Since $X_T^\Delta$ is strongly convergent,
$$\lim_{M \to \infty} E\left|X_T - X_T^{\Delta/k^M}\right| = 0.$$

Using (22), we obtain
$$E|X_T - X_T^\Delta| < \sum_{i=0}^\infty C_2 \frac{\Delta^\gamma}{k^{i\gamma}} = C_2 \Delta^\gamma \sum_{i=0}^\infty \frac{1}{k^{i\gamma}} = C_2 \Delta^\gamma \frac{1}{1-k^{-\gamma}} \quad \text{(since } k^{-\gamma} < 1\text{)}$$
$$= C_1 \Delta^\gamma,$$

where $C_1 = \dfrac{C_2}{1 - k^{-\gamma}}$. This completes the proof. □

The above theorem is particularly useful because it enables the numerical estimation of the strong convergence rate without requiring the analytical solution $X_T$ using the following algorithm.

1. Consider a sequence of step sizes $\Delta_r = \frac{\Delta}{k^{r-1}}, r = 1, 2, \ldots, R$, with corresponding numbers of time steps $N_r = k^{r-1} N$, $r = 1, \ldots, R$.
2. Let $\tilde{E}(\Delta_r)$ denote the Monte Carlo estimate of $E\left|X_T^{\Delta/k^{r-1}} - X_T^{\Delta/k^r}\right|$, $r = 1, \ldots, R$.. This can be obtained via similar steps as defined at the beginning of this subsection, with Step 5. being replaced by $Error(\Delta_r) = \left|X_T^{\Delta/k^{r-1}} - X_T^{\Delta/k^r}\right|$.
3. The strong convergence rate $\gamma$ can be estimated as the slope of the least-squares regression line in $\log \tilde{E}(\Delta_r) = \alpha + \gamma \log(\Delta_r) + \varepsilon_r$.

Thus, this method does not require knowledge of the analytical solution $X_T$.

Next, we state that the above theorem also holds for the weak order of convergence and for convergence in mean squared error.

**Theorem 2.2.** *A time-discrete approximation $X_T^\Delta$ with time step $\Delta$ for the solution $X_T$ of an SDE, constructed using a Wiener process, satisfies the following weak convergence property:*

$$\left|E(X_T) - E(X_T^\Delta)\right| < K\Delta^\beta, \tag{25}$$

*for all $\Delta > 0$, where $K$ is independent of $\Delta$, if and only if the discretization is weakly convergent and, for all $k > 1$, there exists a positive constant $K_1$, independent of $\Delta$, such that*

$$\left|E(X_T^\Delta) - E(X_T^{\Delta/k})\right| < K_1 \Delta^\beta, \tag{26}$$

*where $X_T^{\Delta/k}$ denotes the approximation obtained using time step $\Delta/k$ and the same Wiener process.*

*Proof.* The proof is analogous to that of Theorem 2.1 and is therefore deferred to Appendix A, Section A.1. □

**Theorem 2.3.** *A time-discrete approximation $X_T^\Delta$ with time step $\Delta$ for the solution $X_T$ of an SDE, constructed using a Wiener process, satisfies the following mean squared error convergence:*

$$E\left[(X_T - X_T^\Delta)^2\right] < C_1 \Delta^\gamma, \tag{27}$$

*for all $\Delta > 0$, where $C_1$ is independent of $\Delta$, if and only if the discretization is mean-square convergent and, for all $k > 1$, there exists a positive constant $C_2$, independent of $\Delta$, such that*

$$E\left[(X_T^\Delta - X_T^{\Delta/k})^2\right] < C_2 \Delta^\gamma, \tag{28}$$

*where $X_T^{\Delta/k}$ denotes the approximation using time step $\Delta/k$ and the same Wiener process.*

*Proof.* The proof follows the same arguments as in Theorem 2.1 and is therefore presented in Appendix A, Section A.2. □

## 2.5 Numerical study for rate of convergence

Consider the Black–Scholes SDE [3], which is widely used in the financial world

$$dX_t = rX_t\,dt + \sigma X_t\,dW_t. \tag{29}$$

This SDE models the price of a single risky asset. Here $r$ is the mean rate of return, $\sigma > 0$ is the volatility parameter governing the magnitude of random fluctuations, and $(W_t)$ is a Wiener process. The stochastic differential is interpreted in the Itô sense. The corresponding stochastic integral equation is $X_t = X_0 + \int_0^t rX_s\,ds + \int_0^t \sigma X_s\,dW_s$, where the second integral is the Itô integral.

A closed-form solution is known:

$$X_t = X_0 \exp\left(\sigma W_t + \left(r - \tfrac{1}{2}\sigma^2\right)t\right), \tag{30}$$

and the process $(X_t)$ is a geometric Brownian motion. In particular, $X_t$ is lognormally distributed for every $t > 0$. For a time grid with step size $\Delta = t_{n+1} - t_n$, the Euler–Maruyama (EM) scheme is

$$X_{n+1} = X_n + rX_n\Delta + \sigma X_n \Delta W_n, \tag{31}$$

and the Milstein scheme is

$$X_{n+1} = X_n + rX_n\Delta + \sigma X_n \Delta W_n + \tfrac{1}{2}\sigma^2 X_n(\Delta W_n^2 - \Delta). \tag{32}$$

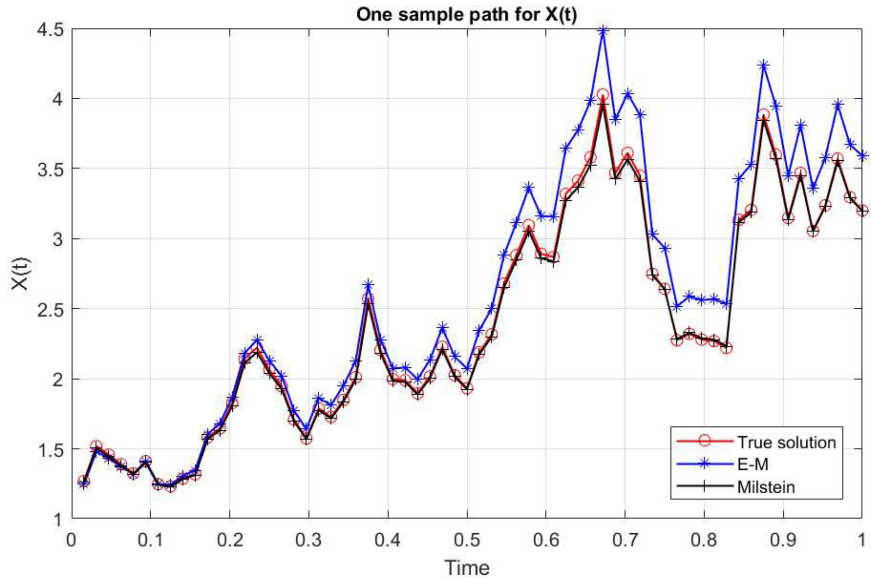

Figure 1: A sample path of the Black–Scholes SDE with $\Delta = 0.0156$ and $N_\Delta = 64$.

In the numerical experiments we set $r = 2$, $\sigma = 1$, and initial value $X_0 = 1$. We use the numerical discretization grid as $\Delta = 0.0313$, i.e., $N_\Delta = 32$. Figure 1 shows one simulated path comparing the Euler and Milstein approximations to the true solution. The Milstein method tracks the analytical solution significantly more closely across all time points.

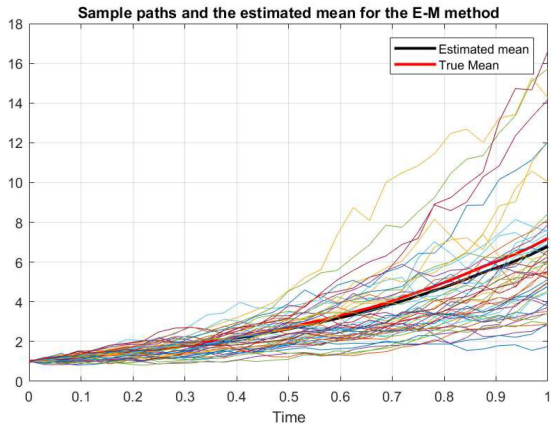
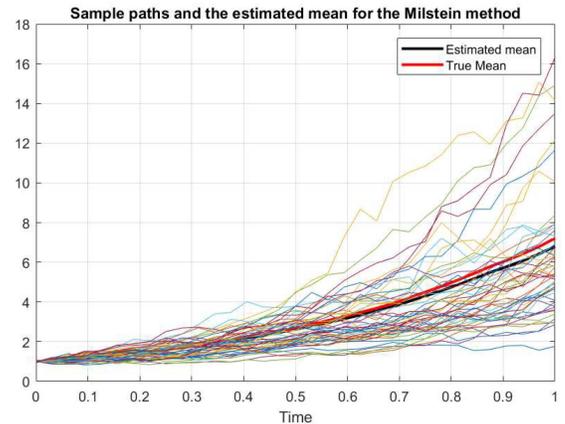

(a) Euler–Maruyama Scheme  (b) Milstein Scheme

Figure 2: Fifty sample paths with true and estimated mean under the Black–Scholes SDE using $\Delta = 0.0313$ ($N_\Delta = 32$).

Figures 2(a) and 2(b) present the empirical mean from 50 Monte Carlo paths using the Euler–Maruyama and Milstein schemes, respectively, compared with the true mean. Both methods perform well for this coarse discretization, though Milstein is slightly more accurate. Figure 3(a) shows the strong convergence rate when the true solution is used as the benchmark. The esti-mated convergence orders are $\gamma_{\text{EM}} = 0.5410, \gamma_{\text{Milstein}} = 0.9596$, which agree closely with the theoretical values 0.5 and 1, respectively. Figure 3(b) uses the approximation defined in Theorem 2.1, comparing $X_T^\Delta$, $X_T^{\Delta/2}$,..., etc. The estimated orders $\gamma_{\text{EM}} = 0.4825, \gamma_{\text{Milstein}} = 0.9441$, again align well with both theory and the estimates obtained using the analytical solution. This also provides a numerical verification of Theorem 2.1.

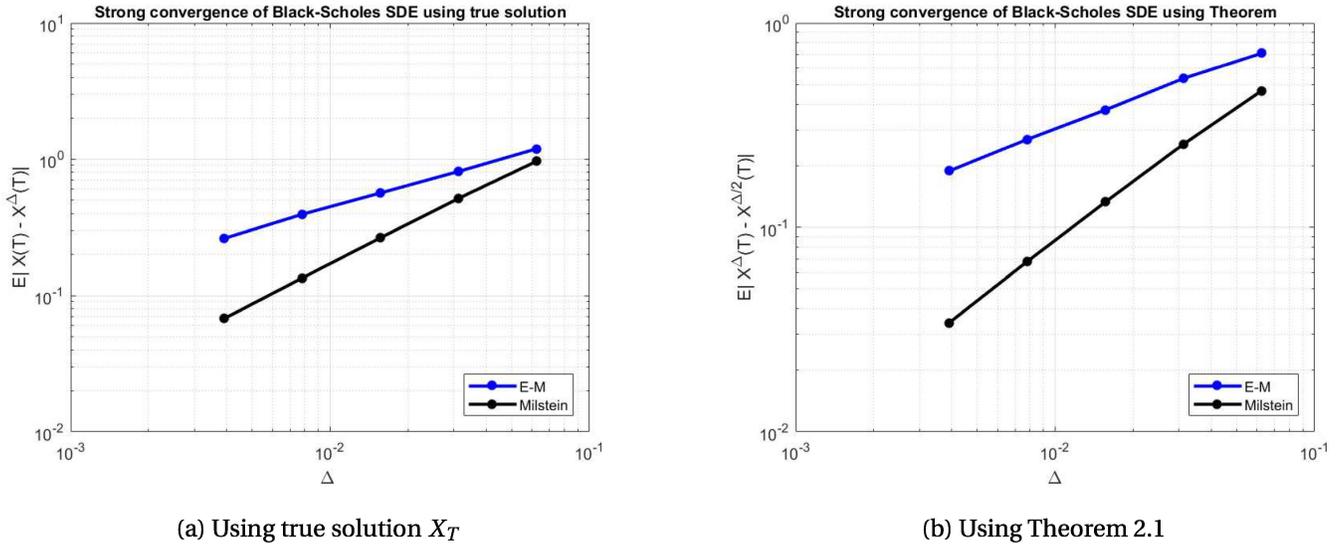

(a) Using true solution $X_T$            (b) Using Theorem 2.1

Figure 3: Strong convergence rates for Euler–Maruyama and Milstein schemes (log–log scale). Expected values estimated using 1000 Monte Carlo paths.

| Method | Using true solution | | Using Theorem 2.1 | |
|---|---|---|---|---|
| | log($C$) | $\gamma$ | log($C$) | $\gamma$ |
| E–M | 1.6729 | 0.5410 | 1.0210 | 0.4825 |
| Milstein | 2.6442 | 0.9596 | 1.8842 | 0.9441 |

Table 1: Estimated strong convergence rate for the Black–Scholes model. Slope $\gamma$ and intercept log($C$) obtained from the least-squares fit of log(Error) vs. log($\Delta$).

## 3 Numerical methods for Multidimensional SDEs

We now extend the numerical framework to multidimensional stochastic differential equations. In many applications, the underlying dynamics are naturally represented by a $d$–dimensional state process $\mathbf{X}_t = (X_t^1, \ldots, X_t^d)^\top$ driven by an $m$–dimensional Wiener process $\mathbf{W}_t = (W_t^1, \ldots, W_t^m)^\top$, where the components $W_t^1, \ldots, W_t^m$ are mutually independent scalar Wiener processes. The multidimensional SDE then takes the general form

$$d\mathbf{X}_t = a(t, \mathbf{X}_t)\, dt + \sum_{j=1}^{m} b^j(t, \mathbf{X}_t)\, dW_t^j, \tag{33}$$

with componentwise representation

$$dX_t^i = a^i(t, \mathbf{X}_t)\, dt + \sum_{j=1}^{m} b^{i,j}(t, \mathbf{X}_t)\, dW_t^j, \quad i = 1, \ldots, d. \tag{34}$$

For notational convenience, this can be written in compact vector form as

$$d\mathbf{X}_t = \sum_{j=0}^{m} b^j(t, \mathbf{X}_t)\, dW_t^j, \tag{35}$$

where we define $b^0(t, x) = a(t, x)$ and set $W_t^0 = t$.

This formulation provides the basis for constructing numerical schemes, such as the Euler–Maruyama and higher-order methods, to approximate trajectories of multidimensional SDE systems. In what follows, we outline these numerical approximations and describe how the extension from the scalar to the vector setting modifies their structure and convergence properties.

## 3.1 Euler–Maruyama Method for Multidimensional SDEs

Since the differential operators associated with a vector-valued SDE act on scalar test functions, the Itô formula is applied componentwise when expanding the solution process $\mathbf{X}(t)$ of (33). Consequently, extending the Euler–Maruyama (E–M) method to multidimensional systems is straightforward: the one-dimensional scheme is applied to each component, with coupling introduced through the drift and diffusion functions.

We first consider the case of a $d$-dimensional state process driven by a scalar Wiener process ($m = 1$). The $i^{\text{th}}$ component of the Euler scheme is then given by

$$X_{n+1}^i = X_n^i + a^i(t_n, \mathbf{X}_n)\Delta_n + b^i(t_n, \mathbf{X}_n)\Delta W_n, \tag{36}$$

for $i = 1, \ldots, d$, where the drift and diffusion coefficients are the vector-valued functions $a = (a^1, \ldots, a^d)$ and $b = (b^1, \ldots, b^d)$. As in the continuous-time system, the components of the discrete-time approximation $\mathbf{X}_n$ remain coupled through these coefficients.

In the general multidimensional setting with $d, m \geq 1$, the $i^{\text{th}}$ component of the Euler–Maruyama update takes the form

$$X_{n+1}^i = X_n^i + a^i(t_n, \mathbf{X}_n)\Delta_n + \sum_{j=1}^m b^{i,j}(t_n, \mathbf{X}_n)\Delta W_n^j, \tag{37}$$

for $i = 1, \ldots, d$. Here, $\Delta W_n^j = W_{t_{n+1}}^j - W_{t_n}^j$ denotes the Wiener increment of the $j^{\text{th}}$ component of the $m$-dimensional standard Brownian motion over the interval $[t_n, t_{n+1}]$. Each increment satisfies $\Delta W_n^j \sim N(0, \Delta_n)$, and increments associated with different components are independent, i.e. $\Delta W_n^{j_1}$ and $\Delta W_n^{j_2}$ are independent whenever $j_1 \neq j_2$. The diffusion coefficient matrix $b = [b^{i,j}]$ is therefore a $d \times m$ matrix, and the Euler–Maruyama method provides a natural and computationally simple discretization for simulating trajectories of multidimensional SDE systems.

## 3.2 Milstein Method for Multidimensional SDEs

To construct higher–order numerical schemes for multidimensional SDEs, we require the differential operators that appear in the multidimensional Itô–Taylor expansion (18). These operators encode the interaction between the drift, diffusion, and noise components and are defined as

$$L^0 = \frac{\partial}{\partial t} + \sum_{k=1}^d a^k \frac{\partial}{\partial x^k} + \frac{1}{2} \sum_{k,l=1}^d \sum_{j=1}^m b^{k,j} b^{l,j} \frac{\partial^2}{\partial x^k \partial x^l},$$

$$L^j = \sum_{k=1}^d b^{k,j} \frac{\partial}{\partial x^k}, \qquad j = 1, \ldots, m, \tag{38}$$

so that each noise component $W_t^j$ corresponds to a distinct operator $L^j$.

With these operators in place, we now generalize the Milstein discretization to SDEs whose state and noise dimensions satisfy $d, m \geq 1$. In this setting, the $i^{\text{th}}$ component of the Milstein update is given by

$$X_{n+1}^i = X_n^i + a^i(t_n, \mathbf{X}_n)\Delta_n + \sum_{j=1}^m b^{i,j}(t_n, \mathbf{X}_n)\Delta W_n^j + \sum_{j_1, j_2=1}^m L^{j_1} b^{i,j_2}(t_n, \mathbf{X}_n) I_{j_1, j_2}[t_n, t_{n+1}], \tag{39}$$

for $i = 1, \ldots, d$. When $j_1 = j_2$, the corresponding iterated integral has the closed-form approximation $I_{j_1, j_1}[t_n, t_{n+1}] = \frac{1}{2}\big((\Delta W_n^{j_1})^2 - \Delta_n\big)$, as given previously in (19).

The primary difficulty in extending the Milstein method beyond the scalar-noise case arises not from the vector-valued state, but rather from the multidimensional Wiener process. When $m \geq 2$, the Milstein update (39) contains mixed stochastic integrals of the form

$$I_{j_1, j_2}[t_n, t_{n+1}] = \int_{t_n}^{t_{n+1}} \int_{t_n}^t dW_s^{j_1} dW_t^{j_2}, \tag{40}$$

for $j_1 \neq j_2$. These mixed terms cannot be written using only the Wiener increments $\Delta W_n^{j_1}$ and $\Delta W_n^{j_2}$, and no simple tractable distribution is available for $I_{(j_1, j_2)}$.

In the following sub-sections, we examine several numerical approaches for approximating these mixed stochastic integrals and analyze how each approach affects the convergence behavior of the multidimensional Milstein scheme.

### 3.2.1 Properties of the Double Itô Integral

We next summarize key properties of the iterated Itô integral (40). For convenience, we denote $I_{(i,j),n} := I_{i,j}[t_n, t_{n+1}]$. As illustrated in Figure 4, each double integral begins at zero and evolves along its own random trajectory. At the end of the time step, the sum of the two mixed integrals satisfies $I_{(i,j),n} + I_{(j,i),n} = \Delta W_n^i \Delta W_n^j$ a.s., $i \neq j$, and their difference defines the stochastic Lévy area, $A_{(i,j),n} := I_{(i,j),n} - I_{(j,i),n}$. The Lévy area plays a central role in stochastic calculus and is fundamental to higher-order numerical schemes.

Using the above identities, each stochastic multiple integral admits the representation

$$I_{(i,j),n} = \tfrac{1}{2}\big(\Delta W_n^i \Delta W_n^j + A_{(i,j),n}\big), \qquad i \neq j. \tag{41}$$

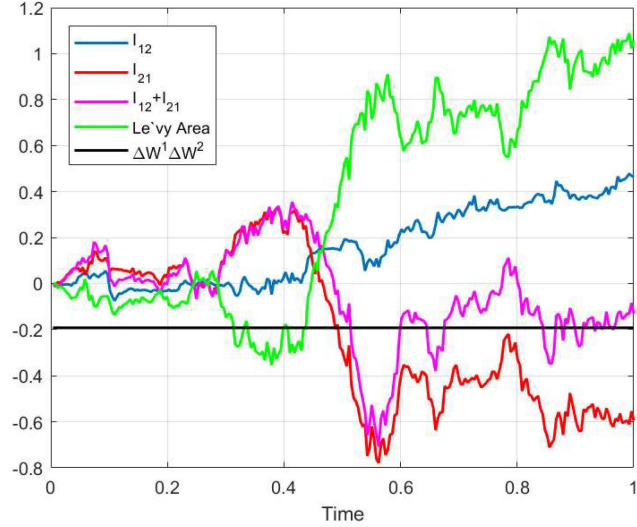

Figure 4: Simulation of a double Itô integral and the associated Lévy area.

The moments of the iterated integrals and the Lévy area are given by

$$E[I_{(i,j),n}] = E[I_{(j,i),n}] = \tfrac{1}{2}\Delta W_n^i \Delta W_n^j, \tag{42}$$

$$\operatorname{Var}(I_{(i,j),n}) = \operatorname{Var}(I_{(j,i),n}) = \tfrac{1}{12}\Delta_n(\Delta_n + R_n^2), \qquad R_n^2 = \Delta W_n^i + \Delta W_n^j, \tag{43}$$

$$E[A_{(i,j),n}] = 0, \qquad \operatorname{Var}(A_{(i,j),n}) = \tfrac{1}{3}\Delta_n(\Delta_n + R_n^2). \tag{44}$$

### 3.2.2 Multidimensional Milstein scheme for special SDE structures

Certain structural properties of an SDE allow the multidimensional Milstein scheme to avoid computing double stochastic integrals (40). Two important cases are *additive noise*, where $b^{i,j}(t,x) \equiv b^{i,j}(t)$, and *commutative noise*, defined by

$$L^{j_1} b^{i,j_2}(t,x) = L^{j_2} b^{i,j_1}(t,x), \qquad 1 \le j_1 \le j_2 \le m,\ i = 1,\ldots,d. \tag{45}$$

For additive noise, all double–integral terms vanish and the Milstein scheme reduces to Euler–Maruyama, yielding strong order $\gamma = 1$. For commutative noise, the identity

$$I_{(j_1,j_2),n} + I_{(j_2,j_1),n} = \Delta W_n^{j_1} \Delta W_n^{j_2}, \qquad j_1 \ne j_2, \tag{46}$$

together with (19), gives the simplified Milstein update

$$X_{n+1}^i = X_n^i + a^i(t_n,X_n)\Delta_n + \sum_{j=1}^m b^{i,j}(t_n,X_n)\Delta W_n^j + \tfrac{1}{2}\sum_{j=1}^m L^j b^{i,j}(t_n,X_n)\bigl((\Delta W_n^j)^2 - \Delta_n\bigr) + \tfrac{1}{2}\sum_{\substack{j_1,j_2=1 \\ j_1 \ne j_2}}^m L^{j_1} b^{i,j_2}(t_n,X_n)\Delta W_n^{j_1}\Delta W_n^{j_2}.$$

A common special case is *diagonal noise* ($d=m$), where each component is driven only by its own Wiener process, i.e. $b^{i,j}(t,x) = 0$ for $i \ne j$ and $\partial b^{i,i}/\partial x^k = 0$ for $k \ne i$. The scheme then decouples componentwise:

$$X_{n+1}^i = X_n^i + a^i(t_n,X_n)\Delta_n + b^{i,i}(t_n,X_n)\Delta W_n^i + \tfrac{1}{2} L^i b^{i,i}(t_n,X_n)\bigl((\Delta W_n^i)^2 - \Delta_n\bigr). \tag{47}$$

Another important structure is *linear noise*, where $b^{i,j}(t,x) = b^{i,j}(t)x^i$, which also satisfies the commutativity condition (45).

## 3.3 Computation of the multiple stochastic integral

In general, the Milstein scheme requires the evaluation of multiple Itô stochastic integrals, except in the special structural cases discussed above. As noted earlier, computing these integrals is challenging, as the task reduces to sampling the triples $(\Delta W_n^i, \Delta W_n^j, A_{(i,j),n})$. While the increments $(\Delta W_n^i, \Delta W_n^j)$ can be sampled directly as independent Gaussian variables, no exact method is known for sampling the Lévy area $A_{(i,j),n}$ conditional on these increments. Several approximation techniques have been proposed in the literature; we review these methods below. In addition, we introduce a new approach based on the solution of an auxiliary multidimensional SDE, and we study its convergence properties and compare its performance with existing methods, both theoretically and numerically.

## 3.4 Lévy–Fourier Method

[13] introduced the Brownian bridge associated with an $m$-dimensional Wiener process, and [11] developed an approximation of multiple stochastic integrals based on a componentwise Fourier expansion of this bridge. Iterative integration with respect to the Wiener components yields the approximation for $I_{(j_1,j_2),n}$:

$$I^p_{(j_1,j_2),n} = \Delta_n \left( \frac{1}{2} \xi_{j_1} \xi_{j_2} + \sqrt{\rho_p}(\mu_{j_1} \xi_{j_2} - \mu_{j_2} \xi_{j_1}) \right) + \frac{\Delta_n}{2\pi} \sum_{r=1}^{p} \frac{1}{r} \left( \zeta_{j_1,r}(\sqrt{2}\xi_{j_2} + \eta_{j_2,r}) - \zeta_{j_2,r}(\sqrt{2}\xi_{j_1} + \eta_{j_1,r}) \right),$$

where $\rho_p = \frac{1}{12} - \frac{1}{2\pi^2}\sum_{r=1}^{p}\frac{1}{r^2}$, $\xi_j = \frac{\Delta W^j_n}{\sqrt{\Delta_n}}$, and $\mu_j, \eta_{j,r}$, and $\zeta_{j,r}$ are mutually independent $N(0,1)$ variables, also independent of $\xi_j$, for $j = 1,\ldots,m$ and $r = 1,\ldots,p$. We refer to this approach as the L–F method. The truncation level $p$ determines the approximation accuracy.

**Result 3.1.** *Assume $\Delta_n = \Delta$ for all $n$. A sufficient condition for strong order $\gamma = 1$ for the parent Milstein scheme is $p \geq K/\Delta$ in the Lévy–Fourier method for some constant $K > 0$.*

*Proof.* [11] established $E\left[(I_{(j_1,j_2)} - I^p_{(j_1,j_2)})^2\right] \leq K_1 \frac{\Delta^2}{p}$. Corollary 10.6.5 of [12] states that strong order 1 is ensured if $E\left[(I_{(j_1,j_2)} - I^p_{(j_1,j_2)})^2\right] \leq K_2 \Delta^3$. Combining the inequalities implies $p > K/\Delta$ with $K = K_1/K_2$. □

The L–F method requires generating $2(2p + 1)$ additional normal random variables for each pair $(j_1, j_2)$, increasing variance and computational cost when $\Delta$ is small. Moreover, because the derivation includes two layers of approximation, accuracy is reliable only for sufficiently small time steps. Consequently, the L–F method is rarely used in practice for computing multiple Itô integrals.

### 3.4.1 E–M Kloeden Subdivsion Method

A simple alternative for approximating the multiple stochastic integral is to embed it into the solution of a suitable auxiliary SDE. Consider the 2-dimensional system

$$dX^1_t = X^2_t \, dW^1_t, \quad dX^2_t = dW^2_t, \tag{48}$$

on the interval $[t_n, t_{n+1}]$. [10] claimed that the solution of (48) with initial values $X^1_{t_n} = 0$, $X^2_{t_n} = W^2_{t_n}$ satisfies

$$X^1_{t_{n+1}} = I_{(2,1)}[t_n, t_{n+1}], \quad X^2_{t_{n+1}} = \Delta W^2_n. \tag{49}$$

Thus, $I_{(2,1),n}$ can be approximated by applying the 2D Euler–Maruyama scheme to (48) over $[t_n, t_{n+1}]$. Let $t'_k = t_n + k\delta_n$ be a subdivision of $[t_n, t_{n+1}]$ with increments $\delta W^j_{n,k} = W^j_{t'_{k+1}} - W^j_{t'_k}$, $k = 0, \ldots, n_K - 1$. The Euler–Maruyama discretization of (48) is

$$Y^1_{k+1} = Y^1_k + Y^2_k \delta W^1_{n,k}, \quad Y^2_{k+1} = Y^2_k + \delta W^2_{n,k}, \tag{50}$$

with initial values $Y^1_0 = 0$ and $Y^2_0 = W^2_{t_n}$. The approximation of the multiple integral is then $I_{(2,1),n} \approx Y^1_{n_K}$. For any pair $j_1 \neq j_2$, the same construction applies by relabeling the components of (48). We refer to this approach as the *E–M Kloeden* method. When used within a Milstein scheme for a parent SDE, this produces a two-level procedure: at each step $[t_n, t_{n+1}]$, one Milstein term is computed via an inner Euler discretization of (48), a strategy often termed a *subdivision method*.

### 3.4.2 E–M IC= 0 Subdivision Method

Although the method proposed in [10] provides an elegant Euler–Maruyama–based approximation of the multiple stochastic integral, it contains an error in the choice of initial conditions. [4] showed that the solution of the SDE (48) satisfies $X^1_{t_{n+1}} = I_{(2,1),n}$, $X^2_{t_{n+1}} = \Delta W^2_n$, when the initial values are $X^1_{t_n} = 0$, $X^2_{t_n} = 0$. Regardless of the initial condition, the increments $\delta W^j_{n,k}$ used in the Euler scheme must be generated via a Brownian bridge such that $\Delta W^j_n = \sum_{k=1}^{n_K} \delta W^j_{n,k}$.

To assess the impact of the initial condition, we performed a numerical study using 5000 sample paths of two independent Wiener processes with step size $\Delta = 0.0625$, each subdivided into 512 finer increments. We considered $\delta_n = \Delta/n_K$ for $n_K = 2^k$, $k = 1, \ldots, 9$. From (46), $I_{(2,1),n} + I_{(1,2),n} = \Delta W^1_n \Delta W^2_n$ a.s., and therefore defining $\text{Error}_n = \Delta W^1_n \Delta W^2_n - (I_{(2,1),n} + I_{(1,2),n})$, the error should converge to zero as $n_K$ increases. Figure 5 presents the Monte Carlo estimates. The original *E–M Kloeden* method converges only at the first interval, while the corrected initial condition $X^1_{t_n} = X^2_{t_n} = 0$ yields convergence at all time points.

Thus, the auxiliary SDE approach of [10] is valid provided the initial condition is replaced by $X^1_{t_n} = 0$, $X^2_{t_n} = 0$. We refer to this corrected scheme as the *E–M IC=0* method. Its accuracy depends on the subdivision parameter $n_K$ (equivalently, on $\delta_n$).

**Result 3.2.** *For the Milstein scheme to achieve strong order $\gamma = 1$, it is sufficient to choose $\delta_n \leq K\Delta^2$, i.e. $n_K > 1/(K\Delta)$, in the E–M IC=0 method for all intervals $[t_n, t_{n+1}]$, where $K > 0$ and $\Delta$ is the Milstein step size.*

*Proof.* Let $I^{EM,\delta}_{(j_1,j_2),n}$ denote the approximation produced by the *E–M IC=0* scheme. From [4], $E\left[(I_{(j_1,j_2),n} - I^{EM,\delta}_{(j_1,j_2),n})^2\right] = \frac{1}{2}\Delta \delta_n$. Corollary 10.6.5 of [12] states that strong order 1 of the Milstein scheme is ensured if, for each interval, $E\left[(I_{(j_1,j_2),n} - I^{EM,\delta}_{(j_1,j_2),n})^2\right] \leq K_1 \Delta^3$. Combining these two results yields $\frac{1}{2}\Delta \delta_n < K_1 \Delta^3$, i.e. $\delta_n < 2K_1 \Delta^2$. Setting $K = 2K_1$ gives the stated result. Since $n_K = \Delta/\delta_n$, the condition is equivalent to $n_K > 1/(K\Delta)$. □

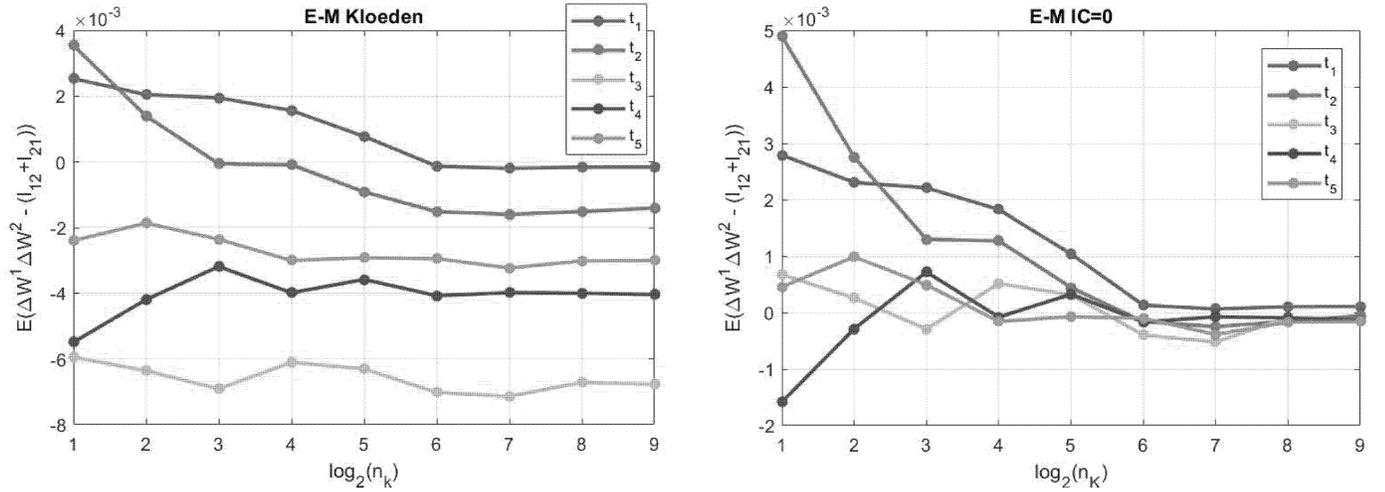

Figure 5: Monte Carlo estimate of the mean error for the first five time intervals, each subdivided into $n_K = 2^k$ parts, $k = 1, \ldots, 9$. Left: *E–M Kloeden* method. Right: *E–M IC=0* method.

### 3.4.3 Approximate Milstein Solution to an SDE

We now propose an alternative approach for approximating the multiple stochastic integral by solving the SDE (48) using an approximate Milstein scheme with initial condition $X^1_{t_n} = 0$, $X^2_{t_n} = 0$. The exact Milstein discretization of (48) on the subintervals $[t_{n,k}, t_{n,k+1}]$ is

$$Y^1_{k+1} = Y^1_k + Y^2_k \delta W^1_{n,k} + I_{(2,1),n,k},$$
$$Y^2_{k+1} = Y^2_k + \delta W^2_{n,k}, \qquad (51)$$

with initial values $Y^1_0 = 0$ and $Y^2_0 = 0$. Since the scheme (51) contains the unknown multiple stochastic integral $I_{(2,1),n,k}$, it cannot be implemented directly. However, from (41), $I_{(i,j),n,k} = \frac{1}{2}\left(\delta W^i_{n,k} \delta W^j_{n,k} + A_{(i,j),n,k}\right)$, $i \neq j$, where $E[A_{(i,j),n,k}] = 0$. If we approximate the Lévy area by $A_{(i,j),n,k} = 0$, the approximate Milstein update becomes

$$Y^1_{k+1} = Y^1_k + Y^2_k \delta W^1_{n,k} + \tfrac{1}{2}\delta W^1_{n,k}\delta W^2_{n,k},$$
$$Y^2_{k+1} = Y^2_k + \delta W^2_{n,k}, \qquad (52)$$

with the same initial condition. We refer to this approximation as the "*Milstein L = 0*" method.

This construction yields a two–stage procedure in which the double Itô integral needed in the parent Milstein scheme is itself approximated by a Milstein-type recursion. As in the "*E–M IC= 0*" approach, the accuracy depends on the choice of $n_k = \Delta/\delta_n$. To preserve strong order $\gamma = 1$ for the outer Milstein method, it again suffices to take $n_k > 1/(K\Delta)$ for some constant $K > 0$. The key advantage of the "*Milstein L = 0*" method is that it attains a smaller mean-square error constant than "*E–M IC= 0*", and hence achieves the same convergence order with roughly half as many subintervals.

**Result 3.3.** *A sufficient condition for achieving strong order $\gamma = 1$ of the parent Milstein scheme when the double integral is approximated using the "Milstein L = 0" method is $\delta_n \leq K\Delta^2$, i.e., $n_k > \frac{1}{K\Delta}$, for all the time interval $[t_n, t_{n+1}]$, $n = 0, 1, 2, \ldots$, where $K$ is a positive constant and $\Delta$ is the time discretization of each step for the parent Milstein scheme.*

*Proof.* Let $I^{M,\delta}_{(j_1,j_2),n}$ denote the approximation obtained from (52). From a result of [4] we have,

$$E\left[\left(I_{(j_1,j_2),n} - I^{M,\delta}_{(j_1,j_2),n}\right)^2\right] = \tfrac{1}{4}\Delta\delta_n. \qquad (53)$$

Corollary 10.6.5 of [12] states that strong order $\gamma = 1$ for the parent Milstein scheme holds provided

$$E\left[\left(I_{(j_1,j_2),n} - I^{M,\delta}_{(j_1,j_2),n}\right)^2\right] \leq K_2\Delta^3. \qquad (54)$$

Combining (53) and (54) yields the sufficient condition $\tfrac{1}{4}\Delta\delta_n < K_2\Delta^3$, i.e. $\delta_n < K\Delta^2$ with $K = 4K_2$, which is equivalent to $n_k > 1/(K\Delta)$, where $n_k = \frac{\Delta}{\delta_n}$. □

**Result 3.4.** *Suppose the "E–M IC= 0" method uses $\delta_n = K\Delta^2$, i.e., $n_k = 1/(K\Delta)$), for all the time interval $[t_n, t_{n+1}]$, $n = 0, 1, 2, \ldots$, to obtain the $\gamma = 1$ order of strong convergence of the Milstein scheme. Then the "Milstein L = 0" method achieves the same strong-order accuracy with $\delta_n = 2K\Delta^2$, i.e. $n_k = \tfrac{1}{2}(1/(K\Delta))$.*

*Proof.* From the results in [4] we have,

$$E\left[(I_{(j_1,j_2),n} - I^{M,\delta}_{(j_1,j_2),n})^2\right] = \tfrac{1}{4}\Delta\delta_n, \tag{55}$$

$$E\left[(I_{(j_1,j_2),n} - I^{EM,\delta}_{(j_1,j_2),n})^2\right] = \tfrac{1}{2}\Delta\delta_n. \tag{56}$$

Also, from Corollary 10.6.5 of [12] it is evident that a sufficient condition to obtain a $\gamma = 1$ strong order convergence of the Milstein scheme is that for each time interval $[t_n, t_{n+1}], n = 0, 1, \ldots$

$$E\left[(I_{(j_1,j_2),n} - I^{\delta}_{(j_1,j_2),n})^2\right] \le K_2 \Delta^3. \tag{57}$$

Thus the "*E–M IC= 0*" method requires $\delta_n < 2K_2\Delta^2$, whereas the "*Milstein L = 0*" method requires only $\delta_n < 4K_2\Delta^2$. Hence, the same Milstein accuracy is obtained by doubling $\delta_n$ (halving $n_k$) for the "*Milstein L = 0*" method. □

The above results show that the "*Milstein L = 0*" method matches the accuracy of "*E–M IC= 0*" while requiring only half as many subintervals, providing an approximately twofold computational speedup. Figure 6 compares the performance of the four approximation methods by plotting the estimated mean absolute error over the last time interval for $k = 1, \ldots, 5$, where $\delta_k = \Delta/2^k$. As expected, the Lévy–Fourier method is insensitive to grid refinement. The *E–M Kloeden* method exhibits large error due to its incorrect initial condition. The *E–M IC= 0* method improves with finer grids, while the proposed *Milstein L = 0* method consistently attains the smallest error across all subdivisions, confirming its superior efficiency.

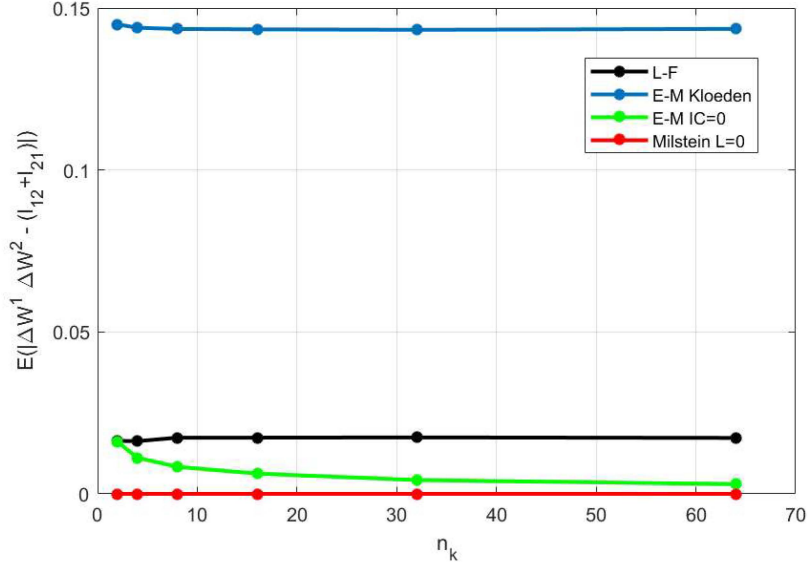

Figure 6: Estimated mean absolute error versus number of subdivision intervals for the four methods. Here $N_\Delta = 32$, and for the "*L–F*" method we take $p = n_k$.

## 4 Applications to Multivariate SDEs in Finance

Having developed efficient numerical methods for multidimensional SDEs, we now illustrate their application to a widely used stochastic volatility model in finance. The classical Black–Scholes model assumes constant volatility and therefore underprices options when volatility fluctuates. Stochastic volatility models address this by allowing the variance itself to evolve randomly. Among these, the Heston model [7] is one of the most influential.

**Generalized Heston Model**: Under the risk-neutral measure, the generalized Heston dynamics are given by the following multidimensional SDE:

$$\begin{aligned} dS_t &= rS_t\,dt + \sqrt{v_t}\,S_t\,d\hat{W}^1_t, \\ dv_t &= \kappa(\theta - v_t)\,dt + \xi v_t^\eta\,d\hat{W}^2_t, \end{aligned} \tag{58}$$

where $S_t$ is the asset price, $v_t$ the instantaneous variance, $r$ the risk-free rate, $\theta$ the long-run mean of $v_t$, $\kappa$ the mean-reversion speed, $\xi$ the volatility-of-volatility, and $\eta$ the elasticity parameter. The Brownian motions $(\hat{W}^1_t, \hat{W}^2_t)$ have correlation $\rho$. To work with independent noises, we apply the standard transformation $d\hat{W}^1_t = dW^1_t, d\hat{W}^2_t = \rho\,dW^1_t + \sqrt{1-\rho^2}\,dW^2_t$, which yields the equivalent representation

$$\begin{aligned} dS_t &= rS_t\,dt + S_t\sqrt{v_t}\,dW^1_t, \\ dv_t &= \kappa(\theta - v_t)\,dt + \rho\xi v_t^\eta\,dW^1_t + \sqrt{1-\rho^2}\,\xi v_t^\eta\,dW^2_t. \end{aligned} \tag{59}$$

For numerical approximation, the Euler scheme becomes

$$S_{n+1} = S_n + S_n r \Delta + S_n \sqrt{v_n} \Delta W_n^1,$$
$$v_{n+1} = v_n + \kappa(\theta - v_n)\Delta + \xi v_n^\eta \left(\rho \Delta W_n^1 + \sqrt{1-\rho^2} \Delta W_n^2\right). \tag{60}$$

Applying the Milstein correction to each equation of 59 while ignoring the stochastic variation of $v_t$ in the diffusion of $S_t$ gives the one–dimensional Milstein form

$$S_{n+1} = S_n + S_n r \Delta + S_n \sqrt{v_n} \Delta W_n^1 + \tfrac{1}{2} S_n v_n \left((\Delta W_n^1)^2 - \Delta\right),$$
$$v_{n+1} = v_n + \kappa(\theta - v_n)\Delta + \xi v_n^\eta \Delta W_n^2 + \tfrac{1}{2} \eta \xi^2 v_n^{2\eta-1}\left((\Delta W_n^2)^2 - \Delta\right). \tag{61}$$

The full multidimensional Milstein scheme, which incorporates the multiple stochastic integral $I_{(2,1),n}$, is

$$S_{n+1} = S_n + S_n r \Delta + S_n \sqrt{v_n} \Delta W_n^1 + \tfrac{1}{2}\sqrt{1-\rho^2}\xi S_n v_n^{\eta-\tfrac{1}{2}} I_{(2,1),n} + \left(\tfrac{1}{2} S_n v_n + \tfrac{1}{4}\rho\xi S_n v_n^{\eta-\tfrac{1}{2}}\right)\left((\Delta W_n^1)^2 - \Delta\right),$$
$$v_{n+1} = v_n + \kappa(\theta - v_n)\Delta + \xi v_n^\eta \left(\rho \Delta W_n^1 + \sqrt{1-\rho^2}\Delta W_n^2\right) + \tfrac{1}{2}\eta\xi^2 v_n^{2\eta-1}\left((\rho \Delta W_n^1 + \sqrt{1-\rho^2}\Delta W_n^2)^2 - \Delta\right).$$

A derivation of (??) is provided in Appendix B. For numerical experiments, we use the parameters $r = 0.04$, $\theta = 0.07$, $\kappa = 3$, $\xi = 0.24$, $\rho = 0.1$, $\eta = 2/3$, and initial conditions $S_0 = 100$, $v_0 = 0.25$.

We now apply the numerical methodology developed in Section 3.3 to numerically solve the generalized Heston stochastic volatility model. This model provides a natural testbed for evaluating the performance of multidimensional Milstein-type schemes, as its dynamics involve nonlinear diffusion, correlated Brownian motions, and practical interest in pathwise accuracy for both state variables and option-based functionals.

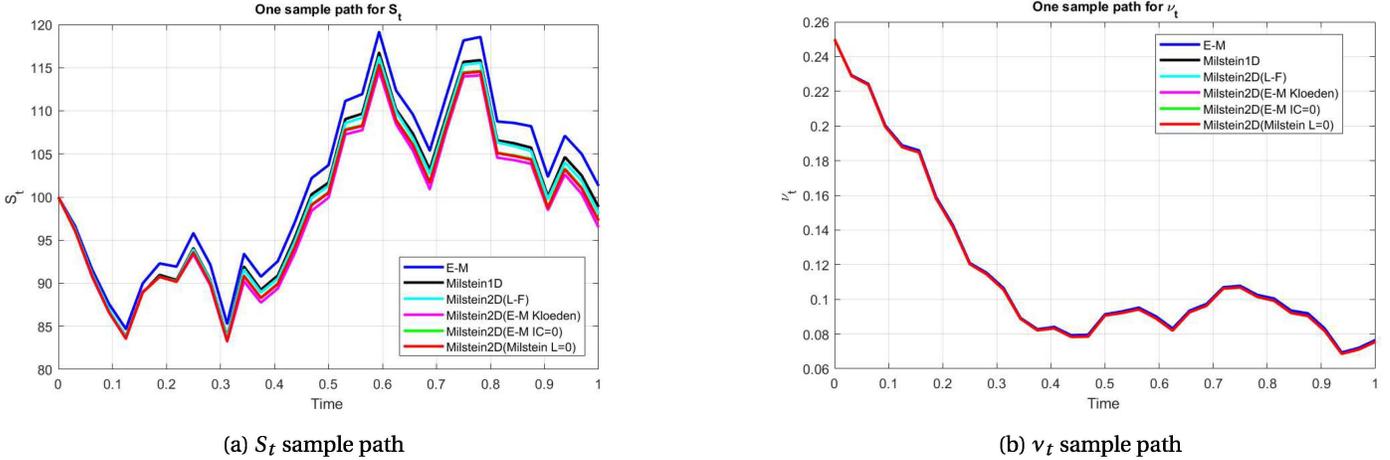

(a) $S_t$ sample path

(b) $v_t$ sample path

Figure 7: Sample paths of $(S_t, v_t)$ generated by different numerical schemes for the generalized Heston model.

Figure 7 (a) shows a representative sample path of $S_t$ for the generalized Heston model computed using the numerical schemes described in Section 3.3. All methods coincide on the initial time grid but gradually diverge as time evolves. Except for Euler–Maruyama, Milstein 1D, and Milstein 2D(L–F), the remaining schemes yield closely aligned trajectories throughout the interval. Figure 7 (b) presents the corresponding sample path for the variance process $v_t$. In contrast to $S_t$, all methods produce nearly indistinguishable approximations, indicating greater numerical stability for the volatility component.

We next compare the strong convergence rates of the stock price $S_t$ and volatility $v_t$ for the numerical schemes introduced in Section 3.3, applied to the generalized Heston model; see Figure 8. Following Theorem 2.1, strong errors are evaluated by comparing terminal values computed with time steps $\Delta$, $\Delta/2$, ..., namely $(S_T^\Delta, v_T^\Delta)$, $(S_T^{\Delta/2}, v_T^{\Delta/2})$, ..., etc. The estimated expected absolute error is plotted against $\Delta$ on a log–log scale, with expectations approximated via Monte Carlo averaging over 1000 sample paths. For all Milstein 2D schemes, the auxiliary subdivision uses $\delta = \Delta^2/2$, except for the Milstein $L = 0$ method, which employs $\delta = \Delta^2$.

For the asset price $S_T$ (Figure 8 (a)), the Milstein 2D(E–M, $IC = 0$) and Milstein 2D(Milstein, $L = 0$) schemes achieve the smallest errors and exhibit slopes close to one. This observation is quantitatively confirmed in Table 2, where these two methods attain the largest estimated convergence rates, $\gamma \approx 1.04$, together with the smallest intercepts $\log(C)$. All remaining schemes display substantially lower slopes, reflecting suboptimal strong convergence. Although the two leading Milstein 2D methods yield comparable accuracy, the Milstein $L = 0$ scheme is approximately twice as fast, since it requires only $\delta = \Delta^2$ (a courser time grid), rather than $\delta = \Delta^2/2$. This computational advantage, achieved without loss of accuracy, confirms the efficiency predicted by Result 3.4. Euler–Maruyama performs worst, with the smallest slope and largest error constant, followed by Milstein 2D(E–M Kloeden), Milstein 2D(L–F), and Milstein 1D.

For the variance process $v_T$ (Figure 8 (b)), Table 2 shows that all methods except Euler–Maruyama achieve strong convergence rates close to one, with nearly identical slopes and intercepts. This behaviour reflects the smoother, mean-reverting structure of the volatility equation and indicates that higher-order corrections primarily benefit the asset price component.

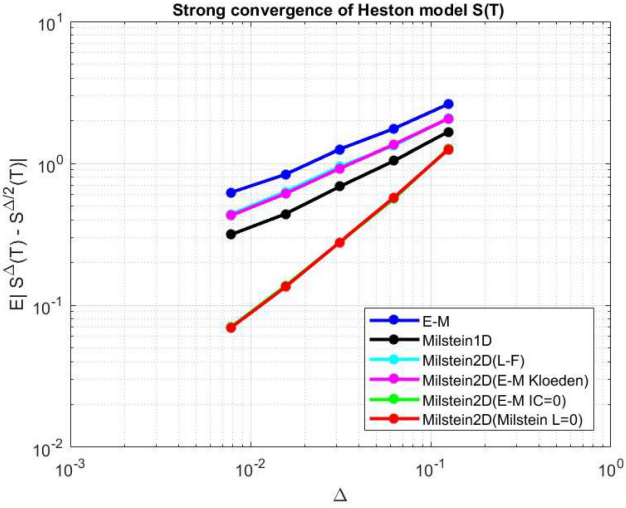
(a) $S(T)$ strong error

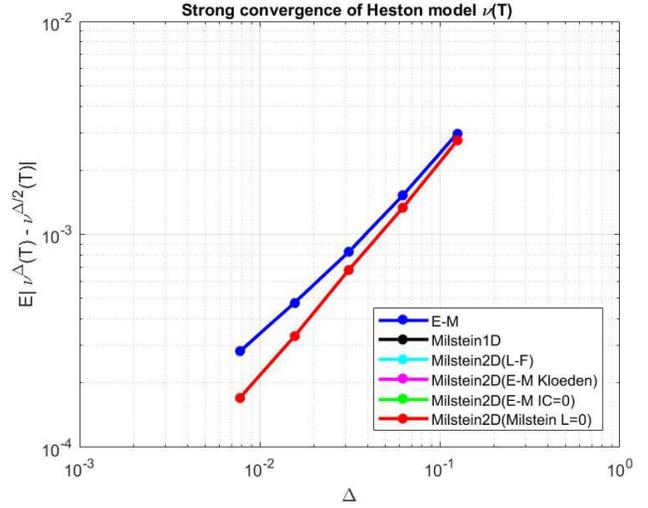
(b) $\nu(T)$ strong error

Figure 8: Strong convergence of $S_T$ and $\nu_T$ under EM, Milstein 1D, and Milstein 2D schemes (log–log scale). For Milstein 2D schemes, $\delta = \Delta^2/2$, except for "Milstein $L = 0$", where $\delta = \Delta^2$.

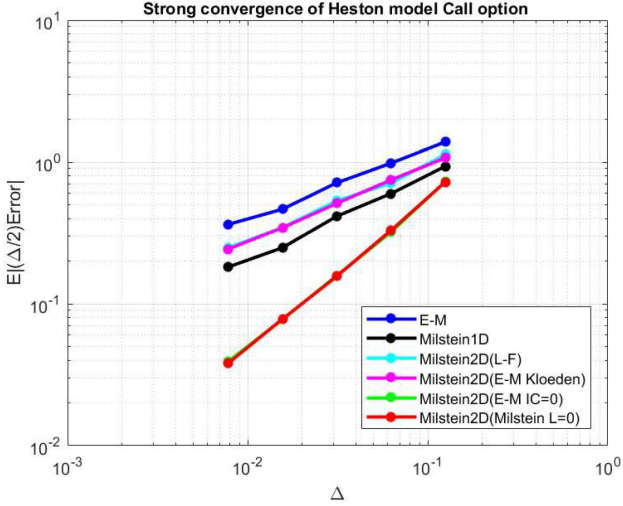
(a) Call option value

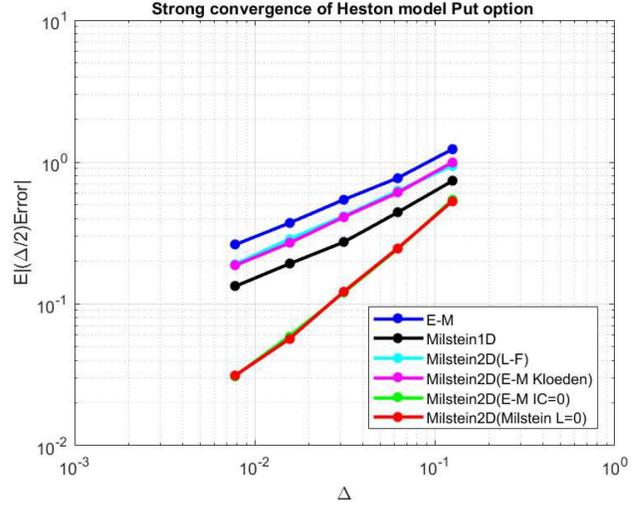
(b) Put option value

Figure 9: Strong convergence of European call and put prices under the generalized Heston model.

Figure 9 extends the comparison to European call and put options. The error curves again mirror those of $S_T$, with Milstein 2D(E–M, $IC = 0$) and Milstein 2D(Milstein, $L = 0$) delivering the smallest errors. Table 2 reports the corresponding slopes and intercepts for option prices and confirms near-first-order strong convergence for these two methods, while all other schemes exhibit significantly lower rates. As in the state-variable case, the Milstein $L = 0$ method achieves this accuracy with half the number of subdivisions, making it the most efficient scheme overall.

Figure 10 reports the strong convergence behaviour of the generalized Heston model when errors are measured in the joint $L_2$ norm for the state vector $X = (S, \nu)$ under expectation. The results are fully consistent with those observed for the individual components and option payoffs in Figures 8 and 9. In particular, the highest convergence rates are again achieved by Milstein 2D(E–M, $IC = 0$) and Milstein 2D(Milstein, $L = 0$), with the latter retaining its computational advantage by using the coarser subdivision $\delta = \Delta^2$ instead of $\delta = \Delta^2/2$. To avoid redundancy, we note that the least-squares fits of the log–log errors yield slopes and intercepts closely matching those reported for $S_T$ in Table 2, confirming the robustness of these conclusions across different error metrics.

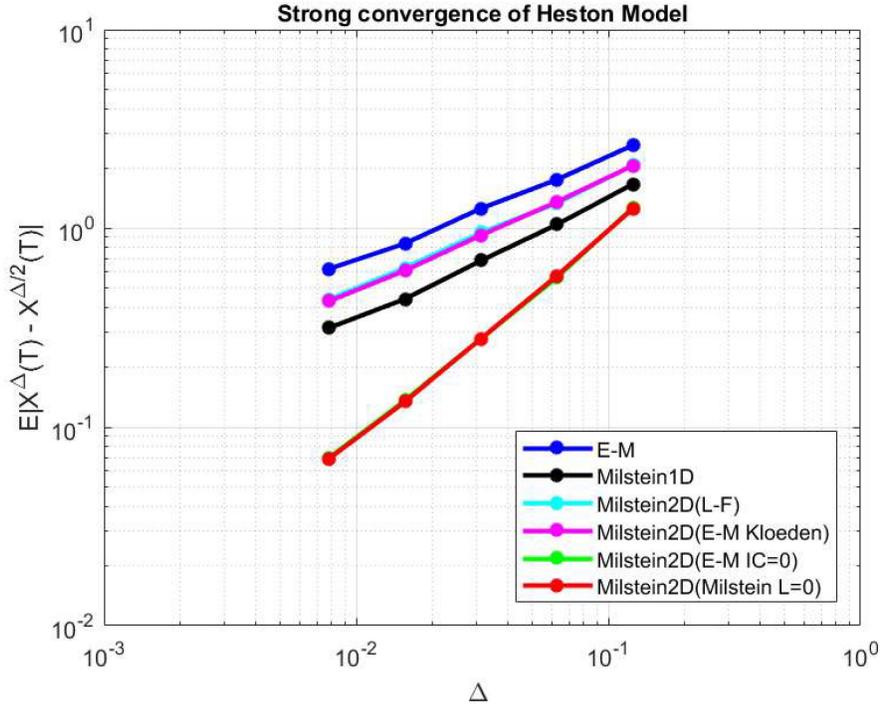

Figure 10: Strong convergence in the $L_2$ norm for $X = (S, v)$ under EM, Milstein 1D, and Milstein 2D schemes. Same subdivision rules as in Figures 8–9.

| Methods | Asset Price $S_T$ | | Variance $v_T$ | | Call Options $V_T$ | | Put Options $V_T$ | |
|---|---|---|---|---|---|---|---|---|
| | log($C$) | $\gamma$ | log($C$) | $\gamma$ | log($C$) | $\gamma$ | log($C$) | $\gamma$ |
| E–M | 2.0246 | 0.5206 | -4.1083 | 0.8475 | 1.3511 | 0.4942 | 1.3196 | 0.5537 |
| Milstein1D | 1.7387 | 0.6058 | -3.8249 | 1.0025 | 1.1647 | 0.5993 | 0.9090 | 0.6143 |
| Milstein2D(*L–F*) | 1.8918 | 0.5685 | -3.8249 | 1.0025 | 1.1966 | 0.5402 | 1.2081 | 0.6022 |
| Milstein2D(*E–M Kloeden*) | 1.8667 | 0.5571 | -3.8249 | 1.0025 | 1.2164 | 0.5419 | 1.1319 | 0.5756 |
| Milstein2D(*E–M IC = 0*) | 2.3461 | 1.0385 | -3.8249 | 1.0025 | 1.8064 | 1.0463 | 1.4727 | 1.0283 |
| Milstein2D(*Milstein L = 0*) | 2.3617 | 1.0441 | -3.8249 | 1.0025 | 1.8391 | 1.0564 | 1.4667 | 1.0282 |

Table 2: Estimated strong convergence rates from least-squares fits of log(Error) versus log($\Delta$). Here $\delta = \frac{\Delta^2}{2}$ i.e. $n_k = 2\Delta^{-1}$ for all the Milstein 2D schemes except "*Milstein L= 0*", for which $\delta = \Delta^2$, i.e., $n_k = \Delta^{-1}$.

## 5 Conclusion

We conclude by emphasizing the significance of Theorem 2.1, which provides a rigorous framework for assessing strong (and weak) convergence of numerical schemes for SDEs. By comparing solutions across different discretization levels it enables accurate estimation of convergence rates without requiring exact solutions, offering both theoretical insight and practical guidance. We have developed and analyzed a practical framework for the numerical approximation of solutions for multidimensional SDEs involving multivariate Itô integrals, focusing on Milstein-type strong schemes. By reformulating multiple stochastic integrals as solutions of auxiliary low-dimensional SDEs, we identified conditions under which subdivision-based Euler approximations are consistent. In particular, we showed that inappropriate initialization may lead to non-convergent approximations and proposed a corrected construction, *E–M IC*= 0, which restores mean-square convergence. We also proposed a simplified yet highly effective variant, the *Milstein L* = 0 scheme, which preserves strong order one convergence while requiring only half the subdivision resolution used by other Milstein 2D constructions. As a result, the *Milstein L* = 0 scheme attains comparable accuracy at a substantially reduced computational cost, making it particularly attractive for large-scale and high-frequency simulations. We established the sufficient condition linking the subdivision step size to the main discretization parameter, which guarantees strong order one convergence of the parent multidimensional Milstein schemes. This provides clear guidance for balancing accuracy and computational cost.

These theoretical findings were illustrated with the generalized Heston stochastic volatility model, a representative coupled nonlinear system driven by correlated noise. Numerical experiments confirmed that Milstein schemes incorporating the proposed stochastic integral approximation achieve the predicted strong convergence rates for both state variables and option-based functionals. Moreover, a simplified variant, *Milstein L* = 0, retains strong order one convergence while reducing computational complexity.

Overall, our results underscore the central role of accurate stochastic integral approximation in multidimensional strong simulation and provide a flexible, computationally efficient methodology applicable to a broad class of SDEs in finance and related fields.

**Acknowledgments:** We acknowledge the late Prof. Wozbor Wozynsky, who provided initial guidance to P. B. on this project until his sudden passing.

# A  Proofs of Theorems

## A.1  Proof of Theorem 2.2

*Proof.* Assume (25) is true for all $\Delta > 0$, then it must be true that for all $k > 1$

$$\left|E(X_T) - E(X_T^{\frac{\Delta}{k}})\right| < C_1 \frac{\Delta^\gamma}{k^\gamma} \tag{62}$$

Now,

$$\begin{aligned}
\left|E(X_T^\Delta) - E(X_T^{\frac{\Delta}{k}})\right| &= \left|E(X_T) - E(X_T^{\frac{\Delta}{k}}) - E(X_T) + E(X_T^\Delta)\right| \\
&= \left|E(X_T) - E(X_T^{\frac{\Delta}{k}})\right| + \left|(EX_T) - E(X_T^\Delta)\right| \\
&< C_1 \frac{\Delta^\gamma}{k^\gamma} + C_1 \Delta^\gamma = C_2 \Delta^\gamma
\end{aligned} \tag{63}$$

where, $C_2 = C_1(1 + \frac{1}{k^\gamma})$, which is independent of $\Delta$. Also

$$\left|E(X_T) - E(X_T^\Delta)\right| < C_1 \Delta^\gamma \implies \lim_{\Delta \to 0} \left|E(X_T) - E(X_T^\Delta)\right| = 0 \tag{64}$$

Hence, one side of the theorem is proved.

Now assume (26) is true for all $\Delta > 0$ and $X_T^\Delta$ is weak convergent.

$$X_T - X_T^\Delta = X_T - X_T^{\frac{\Delta}{k^M}} + \sum_{i=0}^{M-1} \left(X_T^{\frac{\Delta}{k^{i+1}}} - X_T^{\frac{\Delta}{k^i}}\right) \tag{65}$$

Hence,

$$\left|E(X_T) - E(X_T^\Delta)\right| \leq \left|E(X_T) - E(X_T^{\frac{\Delta}{k^M}})\right| + \sum_{i=0}^{M-1} \left|E(X_T^{\frac{\Delta}{k^{i+1}}}) - E(X_T^{\frac{\Delta}{k^i}})\right| \tag{66}$$

Since $X_T^\Delta$ is weakly convergent,

$$\lim_{M \to \infty} \left|E(X_T) - E(X_T^{\frac{\Delta}{k^M}})\right| = 0 \tag{67}$$

Thus from (22) we have,

$$\left|E(X_T) - E(X_T^\Delta)\right| < \sum_{i=0}^\infty C_2 \frac{\Delta^\gamma}{k^{i\gamma}} = C_2 \Delta^\gamma \sum_{i=0}^\infty \frac{1}{k^{i\gamma}} = C_2 \Delta^\gamma \frac{1}{\left(1 - (\frac{1}{k})^\gamma\right)} = C_1 \Delta^\gamma \tag{68}$$

where, $C_1 = C_2 \frac{1}{\left(1-(\frac{1}{k})^\gamma\right)}$. Hence the proof. □

## A.2  Proof of Theorem 2.3

*Proof.* Assume (27) is true for all $\Delta > 0$, then it must be true that for all $k > 1$

$$E\left(X_T - X_T^{\frac{\Delta}{k}}\right)^2 < C_1 \frac{\Delta^\gamma}{k^\gamma} \tag{69}$$

Now,

$$\begin{aligned}
E\left(X_T^\Delta - X_T^{\frac{\Delta}{k}}\right)^2 = E\left(\left(X_T - X_T^{\frac{\Delta}{k}}\right) - \left(X_T - X_T^\Delta\right)\right)^2 = E\left(X_T - X_T^{\frac{\Delta}{k}}\right)^2 + E\left(X_T - X_T^\Delta\right)^2 &< C_1 \frac{\Delta^\gamma}{k^\gamma} + C_1 \Delta^\gamma \\
&= C_2 \Delta^\gamma
\end{aligned} \tag{70}$$

where, $C_2 = C_1(1 + \frac{1}{k^\gamma})$, which is independent of $\Delta$. Also,

$$E\left(X_T - X_T^\Delta\right)^2 < C_1 \Delta^\gamma \implies \lim_{\Delta \to 0} E\left(X_T - X_T^\Delta\right)^2 = 0 \tag{71}$$

Hence, one side of the theorem is proved.

Now assume (28) is true for all $\Delta > 0$ and $X_T^\Delta$ is mean square convergent.

$$X_T - X_T^\Delta = X_T - X_T^{\frac{\Delta}{k^M}} + \sum_{i=0}^{M-1} \left( X_T^{\frac{\Delta}{k^{i+1}}} - X_T^{\frac{\Delta}{k^i}} \right) \tag{72}$$

Hence,

$$E(X_T - X_T^\Delta)^2 \le E(X_T - X_T^{\frac{\Delta}{k^M}})^2 + \sum_{i=0}^{M-1} E(X_T^{\frac{\Delta}{k^{i+1}}} - X_T^{\frac{\Delta}{k^i}})^2 \tag{73}$$

Since $X_T^\Delta$ is mean square convergent, $\lim_{M \to \infty} E(X_T - X_T^{\frac{\Delta}{k^M}})^2 = 0$. Thus, from (22) we have,

$$E(X_T - X_T^\Delta)^2 < \sum_{i=0}^{\infty} C_2 \frac{\Delta^\gamma}{k^{i\gamma}} = C_2 \Delta^\gamma \sum_{i=0}^{\infty} \frac{1}{k^{i\gamma}} = C_2 \Delta^\gamma \frac{1}{\left(1 - (\frac{1}{k})^\gamma\right)} = C_1 \Delta^\gamma \tag{74}$$

where, $C_1 = C_2 \frac{1}{\left(1 - (\frac{1}{k})^\gamma\right)}$. Hence the proof. $\square$

# B  Derivation of the multidimensional Milstein Scheme in the Heston Model (58)

We will derive the Milstein scheme from the transformed generalized Heston model (59) with an independent Wiener process and then transform back to the original Wiener process. From the definition, we get

$$L^1 = b^{11} \frac{\partial}{\partial S} + b^{21} \frac{\partial}{\partial v} = S\sqrt{v} \frac{\partial}{\partial S} + \rho \xi v^\eta \frac{\partial}{\partial v}, \tag{75}$$

and

$$L^2 = b^{12} \frac{\partial}{\partial S} + b^{22} \frac{\partial}{\partial v} = 0 + \sqrt{1-\rho^2} \xi v^\eta \frac{\partial}{\partial v} \tag{76}$$

So,

$$\begin{aligned}
\sum_{j_1,j_2=1}^{2} L^{j_1} b^{1 j_2} I_{(j_1,j_2),n} &= L^1 b^{11} I_{(1,1),n} + L^1 b^{12} I_{(1,2),n} + L^2 b^{11} I_{(2,1),n} + L^2 b^{12} I_{(2,2),n} \\
&= \left( S\sqrt{v} \frac{\partial}{\partial S} + \rho \xi v^\eta \frac{\partial}{\partial v} \right) (S\sqrt{v}) I_{(1,1),n} + 0 \\
&\quad + \sqrt{1-\rho^2} \xi v^\eta \frac{\partial(S\sqrt{v})}{\partial v} I_{(2,1),n} + 0 \\
&= \left( Sv + \frac{1}{2} \rho \xi S v^{\eta - \frac{1}{2}} \right) I_{(1,1),n} + \frac{1}{2} \sqrt{1-\rho^2} \xi S v^{\eta - \frac{1}{2}} I_{(2,1),n}
\end{aligned} \tag{77}$$

$$\begin{aligned}
\sum_{j_1,j_2=1}^{2} L^{j_1} b^{2 j_2} I_{(j_1,j_2),n} &= L^1 b^{21} I_{(1,1),n} + L^1 b^{22} I_{(1,2),n} + L^2 b^{21} I_{(2,1),n} + L^2 b^{22} I_{(2,2),n} \\
&= \left( S\sqrt{v} \frac{\partial}{\partial S} + \rho \xi v^\eta \frac{\partial}{\partial v} \right) (\rho \xi v^\eta) I_{(1,1),n} \\
&\quad + \left( S\sqrt{v} \frac{\partial}{\partial S} + \rho \xi v^\eta \frac{\partial}{\partial v} \right) (\sqrt{1-\rho^2} \xi v^\eta) I_{(1,2),n} \\
&\quad + \sqrt{1-\rho^2} \xi v^\eta \frac{\partial(\rho \xi v^\eta)}{\partial v} I_{(2,1),n} \\
&\quad + \sqrt{1-\rho^2} \xi v^\eta \frac{\partial(\sqrt{1-\rho^2} \xi v^\eta)}{\partial v} I_{(2,1),n} \\
&= \eta \rho^2 \xi^2 v^{2\eta-1} I_{(1,1),n} + \eta \rho \sqrt{1-\rho^2} \xi^2 v^{2\eta-1} I_{(1,2),n} \\
&\quad + \eta \rho \sqrt{1-\rho^2} \xi^2 v^{2\eta-1} I_{(2,1),n} + \eta(1-\rho^2) \xi^2 v^{2\eta-1} I_{(2,2),n}
\end{aligned} \tag{78}$$

Now, we know

$$I_{(j,j),n} = \frac{1}{2} \left( (\Delta W_n^j)^2 - \Delta_n \right), j = 1,2, \quad \text{and} \quad I_{(1,2),n} + I_{(2,1),n} = \Delta W_n^1 \Delta W_n^2. \tag{79}$$

Substituting (79) in (77) and (78) we get

$$\sum_{j_1,j_2=1}^{2} L^{j_1} b^{1 j_2} I_{(j_1,j_2),n} = \left( \frac{1}{2} Sv + \frac{1}{4} \rho \xi S v^{\eta - \frac{1}{2}} \right) \left( (\Delta W_n^1)^2 - \Delta_n \right) + \frac{1}{2} \sqrt{1-\rho^2} \xi S v^{\eta - \frac{1}{2}} I_{(2,1),n} \tag{80}$$

$$\sum_{j_1,j_2=1}^{2} L^{j_1} b^{2j_2} I_{(j_1,j_2),n} = \frac{1}{2}\eta\xi^2 \nu^{2\eta-1}\left(\left(\rho W_n^1 + \sqrt{1-\rho^2}W_n^2\right)^2 - \Delta_n\right) \tag{81}$$

Substituting (80) and (81) in the Milstein scheme (39) we get (**??**).